\documentclass[12pt,reqno]{amsart}
\usepackage{amsmath}
\usepackage{amssymb}
\usepackage{amsthm}
\usepackage{mathrsfs}
\usepackage{latexsym}

\newtheorem{theorem}{Theorem}[section]
\newtheorem{lemma}{Lemma}[section]
\newtheorem{corollary}{Corollary}[section]
\newtheorem{remark}{Remark}[section]
\newtheorem{definition}{Definition}[section]
\newtheorem{proposition}{Proposition}[section]
\newtheorem{example}{Example}[section]
\newtheorem{assumption}{Assumption}[section]
\numberwithin{equation}{section}
\newcommand{\bth}{\begin{theorem}}
\newcommand{\ethe}{\end{theorem}}

\newcommand{\bre}{\begin{remark}}
\newcommand{\ere}{\end{remark}}

\newcommand{\ble}{\begin{lemma}}
\newcommand{\ele}{\end{lemma}}
\newcommand{\bde}{\begin{definition}}
\newcommand{\ede}{\end{definition}}

\newcommand{\bco}{\begin{corollary}}
\newcommand{\eco}{\end{corollary}}

\newcommand{\bpr}{\begin{proposition}}
\newcommand{\epr}{\end{proposition}}

\newcommand{\bexer}{\begin{exercise}}
\newcommand{\eexer}{\end{exercise}}
\newcommand{\breh}{\begin{hint}}
\newcommand{\ereh}{\end{hint}}

\newcommand{\halmos}{\hfill \qed}

\newcommand{\bexam}{\begin{example}}
\newcommand{\eexam}{\end{example}}

\newcommand{\pr} {{\bf Proof.}}

\newcommand{\bfi}{\begin{fig}}
\newcommand{\efi}{\end{fig}}

\newcommand{\beao}{\begin{eqnarray*}}
\newcommand{\eeao}{\end{eqnarray*}\noindent}

\newcommand{\beam}{\begin{eqnarray}}
\newcommand{\eeam}{\end{eqnarray}\noindent}
\newcommand{\E}{\mathbf{E}}
\newcommand{\PP}{\mathbf{P}}
\newcommand{\nto}{n\to\infty}

\newcommand{\xto}{x\to\infty}

\newcommand{\bF}{\overline{F}}

\newcommand{\bV}{\overline{V}}

\newcommand{\bbr}{{\mathbb R}}

\newcommand{\bbb}{{\mathbb B}}

\newcommand{\bbn}{{\mathbb N}}

\newcommand{\vep}{\varepsilon}

\begin{document}
\title[Multivariate strong subexponential distributions]{Multivariate strong subexponential distributions: properties and applications}

\author[C.D. Passalidis ]{Charalampos  D. Passalidis}

\address{Dept. of Statistics and Actuarial-Financial Mathematics,
University of the Aegean,
Karlovassi, GR-83 200 Samos, Greece}
\email{sasd24009@sas.aegean.gr.}

\date{{\small \today}}

\begin{abstract}
In this paper we introduce and study the classes of multivariate strong and strongly subexponential distributions. Some first properties are verified, as for example a type of multivariate analogue of Kesten's inequality, the closure property with respect to convolution, and the conditional closure property with respect to convolution roots. Next, we establish the the single big jump principle for the randomly stopped sums, under the assumption that the random vectors in the summation belong to the class of multivariate strong subexponential distributions. Here the conditions of the counting random variable are weaker in comparison with them in multivariate subexponential class. Further, we establish uniform asymptotic estimates for the precise large deviations in multivariate set up, both for random and non-random sums, when the distribution of the summands belongs to the class of multivariate strongly subexponential distributions. Finally, we provide an application to a non-standard risk model, with independent and identically distributed claim vectors, from the class of multivariate strong subexponential distributions and in the presence of constant interest force. More concretely, the common counting process of the claim vectors constitutes from inter-arrival times, that are independent but not necessarily identically distributed. Under some additional condition, on the 'heavyness' of the counting process tail, we establish a uniform asymptotic estimate for the finite time ruin probability in this model.
\end{abstract}

\maketitle
\vspace{3mm}



\section{Introduction} \label{sec.P.1}

As a matter of fact the heavy tailed distributions in combination with the presence of dependence, represent a hot topic in the frame of applied probability. One way for modeling simultaneously these two concepts, is through the multivariate heavy tailed distributions. Such a way leads, through relatively simple and reasonable assumptions, to reach more general results, in relation with one-dimensional distributions under the presence of specific dependence structures.

From the pioneer paper \cite{haan:resnick:1981} the multivariate regular variation, symbolically $MRV$, attracted the interest of various topics of applied probability, as in extreme value theory, in time-series, in risk theory, in risk management  and random walks, see \cite{hult:lindskog:mikosch:samorodnitsky:2005}, \cite{resnick:2007}, \cite{buraczewski:damek:mikosch:2016}, \cite{konstantinides:li:2016}, \cite{li:2016}, \cite{bhattacharya:palmowksi:zwart:2022}, \cite{li:2022}, \cite{cheng:konstantinides:wang:2024}, among others.

Although the $MRV$, mostly in its standard form, was functioning as antidote for facing the heavy-tailed distributions in combination with the dependence, in fact it remains enough restrictive for practical purposes, especially in actuarial mathematics, where just the regular variation of marginal distributions is too restrictive, see for example in \cite[p. 205, Fig. 7.2]{konstantinides:2018}. Thus, the study of other classes of heavy-tailed distributions in multivariate set up, seems unavoidable.

For the multivariate subexponential distributions we know four different approaches, by \cite{cline:resnick:1992} (via point processes and vague convergence), \cite{omey:2006}, \cite{samorodnitsky:sun:2016} and by \cite{konstantinides:passalidis:2024c}. From the first three approaches, this by \cite{samorodnitsky:sun:2016} looks to be superior, since it has some good enough properties, that are inspired by the corresponding properties of stable random vectors. However, these first three classes show to approach the subexponentiality through a linear approach of multivariate single big jump principle, while the last one, seems to approach the subexponentiality through a non-linear type of multivariate single big jump principle. Therefore, we consider that the last two approaches are functioning in complementary mode, see in \cite[Sec. 5]{konstantinides:passalidis:2024a} for more details.

In this paper we focus in the multivariate subexponentiality $\mathcal{S}_A$ approach, suggested in \cite{samorodnitsky:sun:2016}, see definition in Section 2 below. As in the uni-variate case, in many applications, as for example random walks, randomly stopped sums, queuing theory, risk theory among others, are needed assumptions slightly weaker of subexponentiality, the same happens in the multivariate case as well. More concretely in this paper we are interested in the presence of multivariate single big jump principle in the randomly stopped sums, and asymptotic equivalencies for the precise large deviations in multivariate set up. In order to study these two problems, we shall introduce the multivariate strong and strongly subexponential distributions, that belong to class $\mathcal{S}_A$. 

The paper is organized as follows. In Section 2, after providing the necessary preliminary concepts for the one-dimensional and multivariate heavy-tailed distributions, we present the multivariate classes of strong (and strongly) subexponential distributions. 

In Section 3, we examine some closure properties of these classes, like the closure with respect to convolution, and with respect to convolution roots of random vectors distributions. 

In Section 4, we are interested to randomly stopped sums
\beam \label{eq.P.1.1}
{\bf S}_{\tau}=\sum_{i=1}^{\tau} {\bf X}^{(i)}\,,
\eeam
where $\tau$ represent a (one-dimensional) discrete random variable, and the $\{{\bf X}^{(i)}\,,\;i\in \bbn\}$, represent a sequence of independent identically distributed (i.i.d.) copies of the generic random vector ${\bf X}$. We deal with the single big jump principle for the multivariate randomly stopped sum, namely
\beam \label{eq.P.1.2}
\PP[{\bf S}_{\tau} \in x\,A] \sim \E[\tau]\,\PP[ {\bf X} \in x\,A] \,,
\eeam
as $\xto$, see in Section 2 for definition of the set $A$.

Under the condition that ${\bf X}$ has distribution from the class of multivariate subexponential distributions $\mathcal{S}_A$, was established relation \eqref{eq.P.1.2} under the conditions $\E[\tau]< \infty$ and $\E[(1+\vep)^{\tau}]< \infty$, for some $\vep>0$, see in \cite[Th. 4.2]{konstantinides:passalidis:2024g}. Here, under the restriction of the distribution of ${\bf X}$, to belongs in the class of multivariate strong subexponential distributions, we establish relation \eqref{eq.P.1.2}, under the only requirement that  $\E[\tau]< \infty$ and that the $\tau$ has distribution with lighter tail than $F_A$, see Section 2 for its definition. This result is reduced into corresponding one-dimensional results under a concrete form of set $A$, see for more details discussions in Section 4.

In Section 5,  we establish a vector type precise large deviation principle, for non-random and random sums, which are defined as follows
\beam \label{eq.P.1.3}
{\bf S}_{n} = \sum_{i=1}^n{\bf X}^{(i)}\,, \qquad {\bf S}_{N(t)} = \sum_{i=1}^{N(t)} {\bf X}^{(i)}\,,
\eeam
where $\{N(t)\,,\; t \geq 0\}$, is a counting process (one-dimensional). Under the condition that the $\{{\bf X}^{(i)}\,,\;i\in \bbn\}$ are i.i.d. copies of ${\bf X}$, with multivariate strongly subexponential distribution and independent of the  $\{N(t)\,,\; t \geq 0\}$, which satisfies some common assumptions in one-dimensional precise large deviations, we give the asymptotic equivalencies for the 
\beao
\PP[{\bf S}_{n} \in x\,A] \,,\;\PP[{\bf S}_{N(t)} \in x\,A] \,,
\eeao
where $\nto$ (or, $t\to \infty$, respectively), uniformly for any $x$ greater than a quantity, connected with the inverse function of the insensitivity function of $F_A$, see for the definitions in Section 2. These results represent generalization even in the one-dimensional subcase.

Finally, in Section 6, we present an application in a multivariate non-standard risk model. We provide uniform asymptotic estimations, with respect to time, for the entrance probability of the discounted aggregate claims in some rare-set, and for the ruin probability over a finite time horizon, in a risk model with common counting process for the $d$-lines of business and in the presence of constant interest force. The common counting process is NOT necessarily renewal one, since we assume that the inter-arrival times, are independent but not necessarily identically distributed. These results, are achieved under the condition that the counting process has 'lighter' tail, than the common distribution of the claim-vectors, in each rare-set, which belongs to the class of multivariate strong subexponential distributions. The proof of this result is based on an alternative multivariate version of Kesten inequality, and hence the multivariate strong subexponentiality is necessary in some sense.

\section{Preliminary results} \label{sec.P.2}

In this section we present the main concepts for the heavy-tailed distributions, which will be used later. We also define two new multivariate distribution classes, whose simple properties are studied, and are provided the necessary inclusions with respect to other existing multivariate distribution classes.

\subsection{Notation} \label{sec.P.2.1}

In what follows the asymptotic relations hold as $\xto$, except it is referred otherwise. All the vectors are denoted in bold script and are $d$-dimensional, for some $d\in \bbn$, namely we have ${\bf x}=(x_1,\,\ldots,\,x_d)^{\top}$, where ${\bf x}^{\top}$, is the transpose of ${\bf x}$. All the random vectors, have distributions with support on non-negative quadrant $\bbr_+^d = [0,\,\infty)^d$. The operations between two vectors ${\bf x}^{(1)}$, ${\bf x}^{(2)}$, are defined by components, for example ${\bf x}^{(1)}\pm {\bf x}^{(2)}=( x_{1}^{(1)} \pm x_{1}^{(2)},\,\ldots,\, x_{d}^{(1)}\pm x_{d}^{(2)})^{\top}$, and for some constant $\lambda>0$, the scalar product is defined as $\lambda\,{\bf x}=(\lambda\,x_1,\,\ldots,\,\lambda\,x_d)^{\top}$. Additionally, ${\bf 0} = (0,\,\ldots,\,0)$ denotes the origin of the axes. Furthermore, for two real numbers $x,\,y$ we define $x\wedge y :=\min\{x,\,y\}$, $x\vee y :=\max\{x,\,y\}$, and $\left\lfloor x \right\rfloor$ denotes the integer part of $x$. For any set $\bbb$, we denote by $\bbb^c$ its complement, by $\overline{\bbb}$ its closed hull, and by $\partial \bbb$ its border.

For any distribution $V$, we denote by $\bV = 1-V$ its tail. For two independent random variables $Z_1$, $Z_2$ with distributions $V_1$, $V_2$ reprectively, their convolution is denoted by $V_1*V_2$, and with $V^{n*}$, we denote the $n$-th order convolution power of $V$.

Some asymptotic notation is also needed. For two positive functions $f$, $g$, we write $f(x) \sim c\,g(x)$ with $c \in (0,\,\infty)$ if it holds
\beao
\lim \dfrac{f(x)}{g(x)}= c\,.
\eeao

We write $f(x) = o[g(x)]$, if it holds 
\beao
\lim \dfrac{f(x)}{g(x)}= 0\,,
\eeao
and $f(x) = O[g(x)]$, if it holds 
\beao
\limsup \dfrac{f(x)}{g(x)}< \infty\,.
\eeao
Further, we denote $f(x) \asymp g(x)$, if it holds 
$f(x)= O[g(x)]$ and $g(x)= O[f(x)]$, simultaneously.

Similar notations hold for ${\bf f}$, ${\bf g}$ positive $d$-dimensional functions, where the notation corresponds now to ${\bf f}(x\,\bbb)$, ${\bf g}(x\,\bbb)$, for $\bbb \in \bbr^d$, with ${\bf 0} \notin \bbb$. For example, we write ${\bf f}(x\,\bbb) \sim c\,{\bf g}(x\,\bbb)$ with $c \in (0,\,\infty)$ if it holds
\beao
\lim \dfrac{{\bf f}(x\,\bbb)}{{\bf g}(x\,\bbb)}= c\,.
\eeao
Correspondingly, for the $(d+1)$-variate, positive functions  ${\bf f}^*$ and ${\bf g}^*$, we denote ${\bf f}^*(x\,\bbb;\,y)\sim {\bf g}^*(x\,\bbb;\,y)$, uniformly for any $y \in \Delta$, for some non-empty set $\Delta$ if it holds
\beao
\lim_{\xto} \sup_{y \in \Delta} \left|\dfrac{{\bf f}^*(x\,\bbb;\,y)}{{\bf g}^*(x\,\bbb;\,y)} -1 \right| = 0\,.
\eeao  
Additionally, we write ${\bf f}^*(x\,\bbb;\,y) \lesssim {\bf g}^*(x\,\bbb;\,y)$, or ${\bf g}^*(x\,\bbb;\,y) \gtrsim {\bf f}^*(x\,\bbb;\,y)$, uniformly for any $y \in \Delta$, if
\beao
\limsup_{\xto} \sup_{y\in \Delta} \dfrac{{\bf f}^*(x\,\bbb;\,y)}{{\bf g}^*(x\,\bbb;\,y)} \leq 1\,
\eeao
where $y$ and $\Delta$ will be specified when it appears.

\subsection{Heavy-tailed distributions} \label{sec.P.2.2}

Now, we give the preliminary concepts for the heavy-tailed distributions in one dimension. Let us restrict ourselves only on distributions $V$ with support in $\bbr_+$, although most of the definitions below can easily be extended to supports on the whole $\bbr$. Next, all the distributions in this subsection are such that it holds $\bV(x)>0$, for any $x \in \bbr$.

At first we say that a distribution $V$ has heavy tail, symbolically $V \in \mathcal{K}$, if it holds
\beao
\int_0^{\infty} e^{\vep x}\,V(dx) = \infty\,,
\eeao
for any $\vep>0$. Further, we say that a distribution $V$ has long tail, symbolically $V \in \mathcal{L}$, if it holds
\beao
\lim \dfrac {\bV(x-a)}{\bV(x)} = 1\,,
\eeao
for any (or, equivalently, for some) $a>0$. It is well-known that for any distribution $V \in \mathcal{L}$, there exists a function $h\,:\,\bbr_+ \to (0,\,\infty)$, such that $h(x) \to \infty$, $h(x)=o(x)$ and $\bV(x\pm h(x)) \sim \bV(x)$, that is called insensitivity function of $V$, see in \cite[Lem. 2.19, Sec. 2.8]{foss:korshunov:zachary:2013}. 

The most famous class of heavy-tailed distributions is the class of subexponential distributions, symbolically $\mathcal{S}$, with many applications in branching processes, queuing theory, random walks, risk theory and renewal theory, see in \cite{athreya:ney:1972}, \cite{teugels:1975}, \cite{asmussen:1998}, \cite{borovkov:borovkov:2008}, \cite{li:tang:wu:2010}, \cite{konstantinides:2018}, among others. We say that a distribution $V$ is subexponential, symbolically $V \in \mathcal{S}$, if it holds
\beao
\lim \dfrac {\overline{V^{n*}}(x)}{\bV(x)} = n\,,
\eeao
for any (or, equivalently, for some) integer $n \geq 2$. These first three distribution classes are introduced in \cite{chistyakov:1964}. 
   
The class of dominatedly varying distributions is introduced in \cite{feller:1969}. We say that a distribution $V$ has dominatedly varying tail, symbolically $V \in \mathcal{D}$, if it holds
\beao
\limsup \dfrac {\bV(b\,x)}{\bV(x)} < \infty\,,
\eeao
for any (or, equivalently, for some) $b \in (0,\,1)$. It is well-known that $\mathcal{D} \not\subsetneq \mathcal{S}$, $\mathcal{S} \not\subsetneq \mathcal{D}$ and $\mathcal{D} \cap \mathcal{S} \equiv \mathcal{D} \cap \mathcal{L} \neq \emptyset$, see in \cite{goldie:1978}.

A much smaller but famous heavy-tailed distribution class, is the class of regular variation. We say that a distribution $V$ has regularly varying tail with index $\alpha \in (0,\,\infty)$, symbolically $V \in \mathcal{R}_{-\alpha}$, if it holds
\beao
\lim \dfrac {\bV(t\,x)}{\bV(x)} = t^{-\alpha}\,,
\eeao
for any $t> 0$. We will denote 
\beao
\mathcal{R}:= \bigcup_{0<\alpha <\infty} \mathcal{R}_{-\alpha}\,,
\eeao
for more details about regularly variyng distributions and its applications see in \cite{bingham:goldie:teugels:1987}, \cite{mikosch:samorodnitsky:2000a}, \cite{konstantinides:mikosch:2005}, \cite{man:tang:2024}, \cite{liu:yi:2025}, among others. We know that the following inclusions are true $\mathcal{R}\subsetneq \mathcal{D}\cap \mathcal{L}\subsetneq \mathcal{S}\subsetneq \mathcal{L}\subsetneq \mathcal{K}$, see for example in \cite[Ch. 2]{leipus:siaulys:konstantinides:2023}.

Next, we present two distribution classes, that we generalized in multivariate set up and play crucial role in the present paper. Let $Z$ be a random variable with finite mean
\beao
\mu_V := \E[Z] = \int_0^{\infty} \bV(y)\,dy < \infty\,.
\eeao
We say that $V$, with $\mu_V <\infty$, belongs to the class of strong subexponential distributions, symbolically $V \in \mathcal{S}^*$, if it  holds
\beao
\int_0^{x} \bV(x-y)\,\bV(y)\,dy \sim 2 \mu_V\,\bV(x)\,.
\eeao 
The strong subexponential class $\mathcal{S}^*$ is introduced in \cite{klueppelberg:1988}, where was noticed, as its main advantage, that if $V \in \mathcal{S}^*$, then the integrated tail distribution belongs to class $\mathcal{S}$, but the inverse does not hold in general. Class $\mathcal{S}^*$ has been used in several applications in random walks, randomly stopped sums etc, in order to face problems, in which class  $\mathcal{S}$ fails, see for example in \cite{foss:zachary:2003}, \cite{korshunov:2018}, \cite{foss:korshunov:palmowski:2024}, \cite{palmowski:pojer:thonhauser:2025} among others.

Finally, under the same condition $\mu_V < \infty$, we say that $V$ belongs to the class of strongly subexponential distributions, symbolically $V \in \mathcal{S}_*$, if it  holds
\beao
\lim \dfrac{\overline{V_u^{2*}}(x)}{\bV_u(x)} =2\,.
\eeao 
uniformly for any $u \in [1,\,\infty)$, where
\beao
\bV_u(x) =1 \wedge \int_x^{x+u} \bV(y)\,dy \,,
\eeao
for any $x\in \bbr_+$. Class $\mathcal{S}_*$ is introduced in \cite{korshunov:2002}, and found many applications in random walks and large deviations. In \cite[Lem. 1]{kaas:tang:2003} and \cite[p. 28]{denisov:foss:korshunov:2004} were proved the following inclusions $\mathcal{S}^* \subsetneq \mathcal{S}_* \subsetneq\mathcal{S}$, and further for all distributions $V$, with finite mean $\mu_V < \infty$, these inclusion can be extended to
\beam \label{eq.P.2.8}
\mathcal{R}\subsetneq \mathcal{D}\cap \mathcal{L}\subsetneq \mathcal{S}^* \subsetneq \mathcal{S}_* \subsetneq \mathcal{S} \subsetneq \mathcal{L}\subsetneq \mathcal{K}\,,
\eeam
see in \cite[Th. 3.2(a)]{klueppelberg:1988} for the second inclusion.

\subsection{Multivariate analogues of heavy-tailed distributions} \label{sec.P.2.3}

Now, we provide some known definitions of heavy-tailed distribution classes in multivariate framework, and further we introduce two new multivariate classes, that will be used in the rest of the paper. We present also some examples and counterexamples for these classes, together with some properties as the multivariate version of the Kesten's lemma. We should remind that all the random vectors have distribution support on the non-negative quadrant.

We start with the most popular multivariate heavy-tailed distribution $MRV$. Let ${\bf X}$ be a random vector with distribution $F$. We say that $F$ belongs to (standard) multivariate regularly varying distributions, if there exists a Radon measure $\mu$, non-degenerate to zero, and some one-dimensional distribution $V \in \mathcal{R}_{-\alpha}$, with index $\alpha \in (0,\,\infty)$, such that
\beao
\lim \dfrac 1{\bV(x)}\PP\left[ {\bf X} \in x\,\bbb \right] = \mu(\bbb)\,,
\eeao
for any $\mu$-continuous Borel set $\bbb \in [0,\,\infty]^d \setminus \{{\bf 0}\}$. Symbolically we write $F\in MRV(\alpha,\,V,\,\mu)$, see in \cite{resnick:2007} for more information about $MRV$ class.

For the rest multivariate distribution classes, we restrict ourselves to a concrete subset of $\bbr_+^d$. Namely, we need the following family of sets
\beao
\mathscr{R} = \{ A \subsetneq \bbr^d\;:\; A\; \text{ open, \, increasing,\,} A^c\; \text{ convex},\;{\bf 0} \not\in \overline{A} \}\,,
\eeao
where a set $A$ is called increasing if for any ${\bf x} \in A$ and ${\bf y} \in \bbr_+^d$, it holds ${\bf x}+{\bf y} \in A$. The family of sets $\mathscr{R}$ is introduced in \cite[Sec. 4]{samorodnitsky:sun:2016}, and it was proved that for any $A \in \mathscr{R}$, there exists a set of vectors $I_A \subsetneq \bbr^d$, such that it holds
\beao
A = \{{\bf x} \in \bbr^d\;:\;{\bf p}^T\,{\bf x}> 1\,,\; \exists \;{\bf p} \in I_A \}\,
\eeao
Further more we find in \cite[Lem. 4.5]{samorodnitsky:sun:2016} that for a random vector  ${\bf X}$ with distribution $F$, the random variable 
\beam \label{eq.P.2.12}
Y_A = \sup\{u\;:\;{\bf X} \in u\,A \}\,,
\eeam
has proper distribution $F_A$, whose tail is given by 
\beam \label{eq.P.2.13}
\bF_A(x) = \PP[{\bf X} \in x\,A ]=\PP\left[ \sup_{{\bf p} \in I_A } {\bf p}^{\top}\,{\bf X} > x \right]\,,
\eeam
for any $x \in \bbr_+$. This way, through relation \eqref{eq.P.2.13} was defined the multivariate subexponentiality on $A$, for some fixed $A \in \mathscr{R}$, symbolically $F \in \mathcal{S}_A$, when $F_A \in \mathcal{S}$.

Similarly, in \cite{konstantinides:passalidis:2024g}, were defined classes $\mathcal{D}_A$, $\mathcal{L}_A$, $(\mathcal{D} \cap \mathcal{L})_A$, when $F_A \in \{\mathcal{D},\,\mathcal{L},\, \mathcal{D}\cap \mathcal{L}\}$, respectively. For these classes we adopt the notation
\beao
\mathcal{B}_{\mathscr{R}} := \bigcap_{A \in \mathscr{R}} \mathcal{B}_A\,,
\eeao
for any distribution class $\mathcal{B} \in \{ \mathcal{D}\cap \mathcal{L},\,\mathcal{D},\,\mathcal{S},\, \mathcal{L}\}$. By definition, these classes keep the inclusions from the uni-variate classes and are related with the $MRV$ through \cite[Prop. 2.1]{konstantinides:passalidis:2024g} as follows
\beam \label{eq.P.2.15}
\bigcup_{0<\alpha <\infty}MRV(\alpha,\,V,\,\mu) \subsetneq (\mathcal{D}\cap \mathcal{L})_{\mathscr{R}} \subsetneq \mathcal{S}_{\mathscr{R}} \subsetneq \mathcal{L}_{\mathscr{R}}\,.
\eeam

Now we are ready to introduce two new multivariate distribution classes.

\bde \label{def.P.2.1}
Let $A \in \mathscr{R}$ be some fixed set and ${\bf X}$ be a random vector with distribution $F$. We say that $F$ belongs to class of multivariate strong (or, strongly) subexponential distributions on $A$, symbolically $F \in \mathcal{S}_A^*$ (or $F \in \mathcal{S}_{*,A}$), if $F_A \in \mathcal{S}^*$ (or, $F_A \in \mathcal{S}_*$, respectively).
\ede

\bre \label{rem.P.2.1}
It is clear, that if $F \in \mathcal{B}_A$, with $\mathcal{B} \in \{\mathcal{S}^*,\,\mathcal{S}_*\}$, it means that $\mu_{F_A} < \infty$, for $A \in \mathscr{R}$. Furthermore, from relations \eqref{eq.P.2.8} and \eqref{eq.P.2.15}, under the condition $\mu_{F_A} < \infty$, for all  $A \in \mathscr{R}$, then the following inclusions hold
\beao
\bigcup_{0<\alpha <\infty}MRV(\alpha,\,V,\,\mu) \subsetneq (\mathcal{D}\cap \mathcal{L})_{\mathscr{R}} \subsetneq \mathcal{S}_{\mathscr{R}}^*  \subsetneq \mathcal{S}_{*,\mathscr{R}} \subsetneq \mathcal{S}_{\mathscr{R}} \subsetneq \mathcal{L}_{\mathscr{R}}\,,
\eeao
with these inclusions, still intact, for the classes $\mathcal{B}_A$ instead of $\mathcal{B}_\mathscr{R}$, under the weaker condition $\mu_{F_A} < \infty$, only for this  $A \in \mathscr{R}$, where $\mathcal{B} \in \{ \mathcal{D}\cap \mathcal{L},\,\mathcal{S^{*}},\,\mathcal{S}_*,\,\mathcal{S},\, \mathcal{L}\}$.
\ere
Remark \ref{rem.P.2.1} indicates that all $MRV$ distributions, under the restriction $\mu_{F_A} < \infty$, for all $A\in\mathscr{R}$, belong to $\mathcal{S}_{\mathscr{R}}^*$. We can also find other examples of distributions, belonging to $\mathcal{S}_{\mathscr{R}}^*$, through \cite[Prop. 4.15]{samorodnitsky:sun:2016}, under the extra condition $\mu_{F_A} < \infty$ for all $A \in \mathscr{R}$, because their proposition give some sufficient conditions for membership to class $ (\mathcal{D}\cap \mathcal{L})_{\mathscr{R}}$.

\bre \label{rem.P.2.2}
\cite{samorodnitsky:sun:2016}  was focused to the property of the stability with respect to non-negative, non-degenerate to zero, linear combinations of class $\mathcal{S}_{\mathscr{R}}$. Namely, the weighted sum
\beam \label{eq.P.2.16b}
\sum_{i=1}^d l_i\,X_i\,,
\eeam
to have distribution in $\mathcal{S}$. This is the case for all the classes $\mathcal{B}_\mathscr{R}$ from Remark \ref{rem.P.2.1}. Indeed, the main reason for this, is that the sets of the form
\beam \label{eq.P.2.17}
A_1= \left\{ {\bf x} \;:\;\sum_{i=1}^d l_i\,x_i >c \right\}\,,
\eeam
with $c>0$, $l_1,\,\ldots,\,l_d \geq 0$ and $\sum_{i=1}^d l_i=1$, belong in family of sets $\mathscr{R}$, and the fact that $\mathscr{R}$ is cone, with respect to positive scalar multiplication (namely if $A \in \mathscr{R}$, then it holds $\lambda\,A \in \mathscr{R} $, for any $\lambda>0$). Furthermore, the $I_{A_1}$, can be chosen as $\{(l_1/c,\,\ldots,\,l_d/c)\}$, see in \cite[Rem. 4.4]{samorodnitsky:sun:2016}. Hence, by definition of $Y_{A_1}$, and further of $\mathcal{B}_{A_1}$, recall relations \eqref{eq.P.2.12} and \eqref{eq.P.2.13}, we finally obtain that if $F \in \mathcal{B}_{A_1}$, then the non-negative, non-degenerate to zero linear combinations in \eqref{eq.P.2.16b} follow distribution from class $\mathcal{B}$. Therefore, if $F \in \mathcal{B}_{\mathscr{R}}$, then the sum in \eqref{eq.P.2.16b} follows distribution from class $\mathcal{B}$.
\ere

An essencial question that follows, is wheather it holds the opposite in the previous remark, namely from two marginals from class $\mathcal{B} \in \{\mathcal{S},\,\mathcal{S}^*\}$ we can always construct the class $\mathcal{B}_{\mathscr{R}}$, or at least some $\mathcal{B}_A$? 

\bexam \label{exam.P.2.1}
Let ${\bf X}=(X_1,\,X_2)$ a random vector with distribution $F$, and marginals $F_1, \,F_2 \in \mathcal{B}$, with  $\mathcal{B} \in \{  \mathcal{S},\, \mathcal{S}^*\} $. We assume that the $X_1,\,X_2$ are non-negative, independent random variables, and we choose as $A$ the set
\beao
A= \left\{ {\bf x} \;:\;\dfrac 12\,x_1 + \dfrac 12 x_2 >1 \right\}\,,
\eeao
hence, we find $I_A=(1/2,\,1/2)$, see Remark \ref{rem.P.2.2}, then  it holds
\beao
\bF_A(x) = \PP\left[\dfrac 12\,X_1 + \dfrac 12 X_2  >x \right]\,,
\eeao
for any $x>0$. In order to get $F \in \mathcal{B}_A $, it is enough to hold $F_A  \in \mathcal{B}$, hence it is enough the inclusion $F_1*F_2 \in \mathcal{B}$. However, class $\mathcal{B}$ is NOT closed with respect to convolution, see in \cite{leslie:1989} for $\mathcal{S}$ and in \cite{konstantinides:leipus:siaulys:2023} for $\mathcal{S}^*$. So we can find two distributions such that $F_1*F_2 \notin \mathcal{B}$, hence $F \not\in \mathcal{B}_A$ and further $F  \not\in \mathcal{B}_{\mathscr{R}}$. 
\eexam

We observe that if in Example \ref{exam.P.2.1} we have $F_1*F_2 \in \mathcal{B}$, then follows $F \in \mathcal{B}_A$. For necessary and sufficient conditions for closure with respect to convolution of classes $ \mathcal{S}$, and $\mathcal{S}^*$ we refer to \cite{leipus:siaulys:2020} and \cite{konstantinides:leipus:siaulys:2023} respectively. 

As we said, the inclusions of Remark \ref{rem.P.2.1} are not trivial. Except of the example of \cite[Prop. 4.15]{samorodnitsky:sun:2016}, which belongs to class $(\mathcal{D} \cap \mathcal{L})_{\mathscr{R}}$, there exists a rich spectrum of examples, that escape from the $MRV$ framework. For the set
\beao
A_1'= \left\{ {\bf x} \;:\;\sum_{i=1}^d \dfrac 1d\,x_i >c \right\}\,,
\eeao
for $c>0$, the \cite[Exam. 4.2, 4.3, 4.4]{konstantinides:liu:passalidis:2026}, under the additional condition $\mu_{F_{A_1'}} < \infty$, belong to class $\mathcal{S}_{A_1'}^*$, since they belong to $(\mathcal{D} \cap \mathcal{L})_{A_1'}$. Further, in Example $4.5$ of that paper, if $F_1 \in MDA(\Lambda)\cap \mathcal{S}^*$ (instead of $F_1 \in MDA(\Lambda)\cap \mathcal{S}$), where by $MDA(\Lambda)$ we denote the maximum domain of attraction of Gumbel distribution, it holds again $F \in \mathcal{S}_{A_1'}^*$.

Here, we provide another example for the class $\mathcal{S}_{A_i'}^*$, for $i=1,\,2$, with
\beao
A_2'= \left\{ {\bf x} \;:\;x_i >b\,,\;\exists \;i=1,\,\ldots,\,d \right\} \in \mathscr{R}\,,
\eeao
with $b>0$. It is worth to notice that this example is not confined in class $(\mathcal{D} \cap \mathcal{L})_{A_i'}$, for $i=1,\,2$, and is easy enough to be checked, since the conditions are based only on marginals distributions, and on the dependence structure of the components of ${\bf X}$. The dependence structure we shall use is a general weak dependence called linearly wide quadrant dependence, symbolically $LWQD$, introduced by \cite{qian:geng:wang:2022}.

\bexam \label{exam.P.2.2}
Let ${\bf X}=(X_1,\,\ldots,\,X_d)^{\top}$ be a non-negative random vector with distribution $F$. We assume that the components of ${\bf X}$ follow a common distribution $G \in \mathcal{S}^*$, and that they are $LWQD$, namely there exists a sequence of positive constants $\{g_n\,,\;n \in \bbn\}$ and there is $x_0>0$, such that it holds
\beao
\dfrac{\PP\left[\sum_{j=1}^{n-1} X_j > x \;|\; X_n > y\right]}{\PP\left[\sum_{j=1}^{n-1} X_j > x \right]} \leq g_{n-1}\,,
\eeao
uniformly for $n=2,\,\ldots,\,d$, $x,\,y \in [x_0,\,\infty)$. Then $F \in \mathcal{S}_{A_1'}^*$ and $F \in \mathcal{S}_{A_2'}^*$. 

Let us show first that $F \in \mathcal{S}_{A_1'}^*$. It is enough to show that
\beao
\bF_{A_1'}(x) = \PP\left[ {\bf X} \in x\,A_1' \right]=\PP\left[\sum_{j=1}^d \dfrac cd \, X_j >x \right]\,,
\eeao  
for $x>0$, belongs to $\mathcal{S}^*$. 

By \cite[Th. 1]{qian:geng:wang:2022}, we obtain
\beam \label{eq.P.2.a}
\bF_{A_1'}(x) \sim \sum_{j=1}^d \PP\left[ \dfrac cd \,X_j >x \right]=d\,\overline{G'}(x)\,,
\eeam
where $G'$ denotes the distribution of $\dfrac cd \,X_j$, for $j=1,\,\ldots,\,d$, and consequently $G' \in \mathcal{S}^*$. Hence from relation \eqref{eq.P.2.a}, in combination with \cite[Cor. 3.26]{foss:korshunov:zachary:2013}, we find that $F \in \mathcal{S}_{A_1'}^*$. 

Now, we proceed to show $F \in \mathcal{S}_{A_2'}^*$. For the $A_2'$, we can choose some $I_{A_2'} = \left\{\dfrac {{\bf e}_1}b,\,\ldots,\,\dfrac {{\bf e}_d}b \right\}$, where the ${\bf e}_i$ represents the vector, whose $i$-th component is unit and the rest components are zero (see \cite[Rem. 4.4]{samorodnitsky:sun:2016} for more about $I_{A}$). Therefore, for $x>0$ we obtain
\beao
\bF_{A_2'}(x) =\PP\left[ {\bf X} \in x\,A_2' \right]=\PP\left[ \bigvee_{i=1}^d \dfrac 1b\,X_j >x \right]\,.
\eeao
Thus, again by \cite[Th. 1]{qian:geng:wang:2022}, we obtain
\beao
\bF_{A_2'}(x) \sim \sum_{j=1}^d \PP\left[ \dfrac 1b\,X_j >x  \right]=d\,\overline{G''}(x)\,.
\eeao
where $G''$ denotes the common distribution of the $\dfrac 1b\,X_j$, for $j=1,\,\ldots,\,d$. So, from $G'' \in \mathcal{S}^*$, we find $F \in \mathcal{S}_{A_2'}^*$.  
\eexam

Next, we present some properties of the new distribution classes, related with closure properties with respect to strong (or, weak) tail equivalence, and with a multivariate version of Kesten's inequality. We adopt the notation $F_i(x\,A):= \PP[{\bf X}^{(i)} \in x\,A]$ and $Y_A^{(i)} :=\sup\{u\;:\;{\bf X}^{(i)} \in x\,A\}$, with distributions $F_A^{(i)}$ respectively, for $i \in \bbn$, for any fixed  $A \in \mathscr{R}$.

\bpr \label{pr.P.2.1}
Let  $A \in \mathscr{R}$ be some fixed set. We consider the random vectors ${\bf X}^{(1)}$, ${\bf X}^{(2)}$ with distributions $F_1$ and $F_2$.
\begin{enumerate}
\item[(i)]
If $F_1 \in \mathcal{S}_A^*$ ( or, $F_1 \in \mathcal{S}_{*,A}$), $F_2 \in \mathcal{L}_A$ and it holds
\beam \label{eq.P.2.18}
F_1(x\,A) \asymp F_2(x\,A)\,,
\eeam
then $F_2 \in \mathcal{S}_A^*$ ( or, $F_2 \in \mathcal{S}_{*,A}$ respectively).
\item[(ii)]
If $F_1 \in \mathcal{S}_A^*$  and it holds
\beam \label{eq.P.2.19}
\lim \dfrac {F_1(x\,A)}{F_2(x\,A)} =c\,,
\eeam
for some $c \in (0,\,\infty)$, then $F_2 \in \mathcal{S}_A^*$.
\end{enumerate}
\epr

\pr~
\begin{enumerate}
\item[(i)]
Since  $F_A^{(1)} \in \mathcal{S}^*$ ( or $F_A^{(1)}  \in \mathcal{S}_{*}$), $F_A^{(2)} \in \mathcal{L}$ and by relation \eqref{eq.P.2.18} is  implied 
that $\bF_A^{(1)}(x) \asymp \bF_A^{(2)}(x) $, we obtain $F_A^{(2)} \in \mathcal{S}^*$ ( or $F_A^{(2)}  \in \mathcal{S}_{*}$), by \cite[Th. 2.1(b)]{klueppelberg:1988} (or, \cite[Lem. 3]{korshunov:2002}, respectively).
\item[(ii)]
Since  $F_A^{(1)} \in \mathcal{S}^*$ and by relation \eqref{eq.P.2.19} follows $\bF_A^{(1)}(x) \sim c\,\bF_A^{(2)}(x) $, we obtain $F_A^{(2)} \in \mathcal{S}^*$ by \cite[Cor. 3.26]{foss:korshunov:zachary:2013}.
\end{enumerate}
~\halmos

Kesten's inequality plays a crucial role in one dimensional set up, and facilitate several applications. In \cite{samorodnitsky:sun:2016}, we  find a multivariate analogue of this inequality, for the class $\mathcal{S}_A$. Here we give an alternative version of multivariate Kesten's inequality for $\mathcal{S}_A^*$, having in mind the the one dimensional analogue which given in \cite[Th. 2]{denisov:foss:korshunov:2010}.

\bpr \label{pr.P.2.2}
Let $A \in \mathscr{R}$ be some fixed set. If ${\bf X},\,{\bf X}^{(1)},\,\ldots,\,{\bf X}^{(n)}$ are i.i.d. random vectors with common distribution $F \in  \mathcal{S}_A^*$, then for any $c> \mu_{F_A}$, there exists some constant $K>0$, such that it holds
\beam \label{eq.P.2.20}
\dfrac {\PP\left[{\bf X}^{(1)}+\cdots +{\bf X}^{(n)} \in x\,A \right]}{\PP\left[{\bf X} \in x\,A \right]}\leq \dfrac K{\PP\left[{\bf X} \in c\,n\,A \right]}\,,
\eeam
for all $n\in \bbn$ and all $x>0$.
\epr

\pr~
Since ${\bf X}$ is non-negative, random vector, we find $\mu_{F_A} =\E[Y_A] \in (0,\,\infty)$. Hence by \cite[Lem. 4.9]{samorodnitsky:sun:2016} and by \cite[Th. 2, appl.]{denisov:foss:korshunov:2010} we obtain 
\beao
\PP\left[{\bf X}^{(1)}+\cdots +{\bf X}^{(n)} \in x\,A \right] &\leq& \PP\left[Y_A^{(1)}+\cdots +Y_A^{(n)} > x \right] \\[2mm]
&\leq& \dfrac {K\,\PP\left[Y_A^{(1)} > x \right] }{\PP\left[Y_A^{(1)} > c\,n \right] }=\dfrac {K\,\PP\left[{\bf X} \in x\,A \right]}{\PP\left[{\bf X} \in c\,n\,A \right]}\,,
\eeao
for all $x>0$, and $n\in \bbn$.
~\halmos

\section{Convolution and convolution roots properties} \label{sec.P.3}

Here we put more attention on closure properties with respect to convolution of random vectors and with respect to convolution roots of random vectors. At first we consider the closure property with respect to convolution.

Let us start with independent random vectors ${\bf X}^{(1)},\,{\bf X}^{(2)}$ with distributions $F_1$ and $F_2$ respectively. We define their convolution, over the set $x\,A$, as follows
\beam \label{eq.P.3.21}
F_1*F_2(x\,A) = \PP\left[{\bf X}^{(1)}+{\bf X}^{(2)}\in x\,A \right]\,.
\eeam
Hence, for some fixed  $A \in \mathscr{R}$, we say that $F_1*F_2 \in \mathcal{B}_A $, if it holds $F_A^*  \in \mathcal{B}$, where the distribution $F_A^*$ is defined as the distribution of
\beam \label{eq.P.3.22}
Y_A^* := \sup \{ u\;:\; {\bf X}^{(1)}+{\bf X}^{(2)}\in u\,A \}\,.
\eeam  

In \cite[Th. 3.4]{konstantinides:passalidis:2024g} we find that if $F_1,\,F_2 \in \mathcal{L}_A$, then $F_A^{(1)}*F_A^{(2)}  \in \mathcal{S}$ if and only if $F_1*F_2 \in \mathcal{S}_A $, namely the closure issue of multivariate convolution is reduced to the one-dimensional one, when the marginal distributions belong to $\mathcal{L}_A$. A similar result is true for the classes $\mathcal{S}_A^*,\,\mathcal{S}_{*,A}$, as we can see below.

\bth \label{th.P.3.1}
Let  $A \in \mathscr{R}$ be some fixed set, and  $F_1,\,F_2 \in \mathcal{L}_A$.
\begin{enumerate}
\item[(i)]
$F_A^{(1)}*F_A^{(2)}  \in \mathcal{S}^*$ if and only if $F_1*F_2 \in \mathcal{S}_A^* $.
\item[(ii)]
$F_A^{(1)}*F_A^{(2)}  \in \mathcal{S}_*$ if and only if $F_1*F_2 \in \mathcal{S}_{*,A}$.
\end{enumerate}
\ethe 

\pr~
We follow the same methodology for both parts, making only corresponding clarifications when needed. Let $\mathcal{B} \in \{\mathcal{S}^*,\,\mathcal{S}_{*}\}$.

$(\Longrightarrow)$. We suppose that $F_A^{(1)}*F_A^{(2)} \in \mathcal{B}$. So we have to show the inclusion $F_A^* \in \mathcal{B}$ (recall also the relations \eqref{eq.P.3.21} and \eqref{eq.P.3.22}). Since $F_A^{(1)}*F_A^{(2)} \in \mathcal{B} \subsetneq \mathcal{S}$, and $F_A^{(1)},\,F_A^{(2)} \in \mathcal{L} $, by \cite[Th. 1.1]{leipus:siaulys:2020} is implied the max-sum equivalence
\beam \label{eq.P.3.23}
\overline{F_A^{(1)}*F_A^{(2)}} (x) \sim \bF_A^{(1)}(x)+\bF_A^{(2)}(x)\,.
\eeam 
From the other side, since ${\bf X}^{(1)},\,{\bf X}^{(2)}$ are non-negative random vectors, and the set $A$ is increasing, we obtain through Bonferroni inequality
\beam \label{eq.P.3.24} \notag
\bF_A^*(x) &=&\PP[Y_A^* > x] = \PP\left[ {\bf X}^{(1)}+{\bf X}^{(2)} \in x\,A\right] \geq \PP\left[\bigcup_{i=1}^2 \{{\bf X}^{(i)} \in x\,A\} \right]\\[2mm] \notag
&\geq& \PP\left[{\bf X}^{(1)} \in x\,A\right]+ \PP\left[{\bf X}^{(2)} \in x\,A \right]- \PP\left[{\bf X}^{(1)} \in x\,A,\;{\bf X}^{(2)} \in x\,A\right]  \\[2mm] \notag
&=& \PP\left[Y_A^{(1)}> x\right]+ \PP\left[Y_A^{(2)} > x \right]- \PP\left[Y_A^{(1)} > x\right]\, \PP\left[Y_A^{(2)} > x\right]\\[2mm]
&\sim& \PP\left[Y_A^{(1)}> x\right]+ \PP\left[Y_A^{(2)} > x \right]\,.
\eeam 
So, from relations \eqref{eq.P.3.23} and \eqref{eq.P.3.24} we find
\beam \label{eq.P.3.25} 
\limsup \dfrac{\overline{F_A^{(1)}*F_A^{(2)}(x)} }{\bF_A^*(x)} \leq \limsup \dfrac{\PP\left[Y_A^{(1)}> x\right]+ \PP\left[Y_A^{(2)} > x \right]}{\PP\left[Y_A^{(1)}> x\right]+ \PP\left[Y_A^{(2)} > x \right]} =1\,.
\eeam
From the other side, by \cite[Prop. 2.4]{konstantinides:passalidis:2024g}, for all $x>0$ we obtain
\beam \label{eq.P.3.26} 
\dfrac{\bF_A^*(x)}{\overline{F_A^{(1)}*F_A^{(2)}}(x)} \leq 1\,.
\eeam
From relations \eqref{eq.P.3.25} and \eqref{eq.P.3.26}, we conclude 
\beam \label{eq.P.3.27} 
\bF_A^*(x)\sim \overline{F_A^{(1)}*F_A^{(2)}}(x)\,.
\eeam
Next, for case $\mathcal{B} \equiv \mathcal{S}^*$, by \eqref{eq.P.3.27}, since $F_A^{(1)}*F_A^{(2)}  \in \mathcal{S}^*$, it follows $F_A^* \in \mathcal{S}^*$, from \cite[Cor. 3.26]{foss:korshunov:zachary:2013}. Therefore, $F_1*F_2 \in \mathcal{S}_A^* $.

For case $\mathcal{B} \equiv \mathcal{S}_*$,  since $F_A^{(1)}*F_A^{(2)}  \in \mathcal{S}_* \subsetneq \mathcal{S}$, by \eqref{eq.P.3.27}, it follows $F_A^* \in \mathcal{S}$, from the closure property with respect to strong tail equivalence of $\mathcal{S}$, see \cite[Th. 3]{teugels:1975}. Therefore, from \eqref{eq.P.3.27}  $F_A^* \in \mathcal{S}$ and $F_A^{(1)}*F_A^{(2)}  \in \mathcal{S}_*$, via \cite[Lem. 3]{korshunov:2002}, we find $F_A^* \in \mathcal{S}_*$. Thus we find $F_1*F_2 \in \mathcal{S}_{*,A}$ 

$(\Longleftarrow)$. Let suppose that $F_1*F_2 \in \mathcal{B}_A$, with $\mathcal{B} \in \{\mathcal{S}^*,\,\mathcal{S}_{*}\}$. Then for ${Y_A^*}'$ independent and identical copy of ${Y_A^*}$, we obtain that since ${\bf X}^{(1)},\,{\bf X}^{(2)}$ are non-negative, it holds
\beam \label{eq.P.3.28}
\limsup \dfrac{\overline{F_A^{(1)}*F_A^{(2)}}(x)}{\bF_A^*(x)} = \limsup \dfrac{\PP\left[Y_A^{(1)}+ Y_A^{(2)} > x \right]}{\PP\left[Y_A^*> x\right]} \leq  \limsup \dfrac{\PP\left[Y_A^{*}+{Y_A^*}' > x \right]}{\PP\left[Y_A^*> x\right]} =2\,, 
\eeam
where in the last step we use the inclusion $\mathcal{B} \subsetneq \mathcal{S}$. Hence, by relations \eqref{eq.P.3.26} and \eqref{eq.P.3.28}, we find
\beam \label{eq.P.3.29}
\overline{F_A^{(1)}*F_A^{(2)}} \asymp \bF_A^*(x)\,. 
\eeam
From $F_A^{(1)},\,F_A^{(2)} \in \mathcal{L}$ it follows $F_A^{(1)}*F_A^{(2)}\in \mathcal{L}$, see in \cite[Th. 3(b)]{embrechts:goldie:1980}. Hence since $F_A^* \in \mathcal{B}$, $F_A^{(1)}*F_A^{(2)}   \in \mathcal{L}$ and relation \eqref{eq.P.3.29} holds, we obtain $F_A^{(1)}*F_A^{(2)}\in \mathcal{B}$, by \cite[Th. 2.1(b)]{klueppelberg:1988} for class $\mathcal{S}^*$, and by \cite[Lem. 3]{korshunov:2002} for class $\mathcal{S}_{*}$. 
~\halmos

We note that in Theorem \ref{th.P.3.1} can be reduced to whole class of sets $\mathscr{R}$ under the next formulation: $F_A^{(1)}*F_A^{(2)}  \in \mathcal{B}$, for all $A \in \mathscr{R}$ if and only if $F_1*F_2 \in \mathcal{B}_{\mathscr{R}} $. Furthermore, Theorem \ref{th.P.3.1} assers that in $\mathcal{L}_A$ the closure property with respect to convolution for classes $\mathcal{S}_A^*,\,\mathcal{S}_{*,A}$ is reduced to one-dimensional convolution closure properties. For further convolution properties of $\mathcal{S^{*}}$, see in \cite[Sec. 3.11]{leipus:siaulys:konstantinides:2023}. However, we do not know about the convolution closure properties of $\mathcal{S}_{*}$. In \cite[Th. 2]{konstantinides:leipus:siaulys:2023} was shown that class $\mathcal{S}^{*}$ does not satisfy the closure property with respect to convolution, from where it follows the next result.

\bco \label{cor.P.3.1}
Let $A \in \mathscr{R}$. Class $\mathcal{S}_A^{*}$ is not closed with respect to convolution. Namely, we can find two distributions $F_1,\,F_2 \in \mathcal{S}_A^{*}$, such that it holds $F_1*F_2 \not\in \mathcal{S}_A^{*}$.
\eco

The next corollary discloses that class $\mathcal{S}_A^{*}$ is closed with respect to convolution power. Although it is elementary conclusion in comparison with the convolution, can be useful later use. Let remind that $F^{n*}(x\,A) = \PP\left[{\bf X}^{(1)}+\cdots +{\bf X}^{(n)}\in x\,A \right]$, when the random vectors ${\bf X}^{(1)},\,\cdots ,\,{\bf X}^{(n)}$ are i.i.d. with common distribution $F$.

\bco \label{cor.P.3.2}
Let $A \in \mathscr{R}$. If $F \in \mathcal{S}_A^{*}$, then for any $n \in \bbn$ it holds $F^{n*} \in \mathcal{S}_A^{*}$.
\eco

\pr~
Since $F \in \mathcal{S}_A^{*} \subsetneq \mathcal{S}_A$, by \cite[Cor. 4.10]{samorodnitsky:sun:2016} it holds $F^{n*}(x\,A) \sim n\,F(x\,A)$, hence, from Proposition \ref{pr.P.2.1}(ii), we obtain $F^{n*}\in \mathcal{S}_A^{*}$.
~\halmos

Now we consider the closure property of the two new classes with respect to convolution roots. For some fixed set $A \in \mathscr{R}$, we say that the class $\mathcal{B}_A$ is closed with respect to convolution roots, if from inclusion  $F^{n*} \in \mathcal{B}_A$, for some integer $n \geq 2$, it is implied that  $F \in \mathcal{B}_A$. If instead of $A$ we take the whole class of sets $\mathscr{R}$, then from $F^{n*} \in \mathcal{B}_A$, for some integer $n \geq 2$ and any $A \in \mathscr{R}$, the  conclusion is  $F \in \mathcal{B}_{\mathscr{R}}$. It is also easy to see that the result of Proposition \ref{pr.P.3.1}, below, holds also in the case of classes $\mathcal{B}_{\mathscr{R}}$, instead of $\mathcal{B}_A$.

In one-dimensional set up we find relative results in \cite{embrechts:goldie:veraverbeke:1979}, \cite{watanabe:yamamuro:2010}, \cite{xu:foss:wang:2015}, \cite{watanabe:2019}, \cite{cui:wang:xu:2022}, \cite{xu:yu:wang:cheng:2024} among others. In \cite[Th. 3.2]{konstantinides:passalidis:2024h}, we find that class $\mathcal{S}_A$ has the closure property with respect to convolution roots, under the condition that $F \in \mathcal{L}_A$. 

\bpr \label{pr.P.3.1}
Let $A \in \mathscr{R}$ be some fixed set.
\begin{enumerate}
\item[(i)]
If $F^{n*} \in \mathcal{S}_A^*$, for some  integer $n \geq 2$, and $F \in \mathcal{L}_A$, then $F \in \mathcal{S}_A^*$.
\item[(ii)]
If $F^{n*} \in \mathcal{S}_{*,A}$, for some  integer $n \geq 2$, and $F \in \mathcal{L}_A$, then $F \in \mathcal{S}_{*,A}$.
\end{enumerate} 
\epr

\pr~
Let consider $\mathcal{B} \in \{\mathcal{S}^*,\,\mathcal{S}_{*}\}$. Since $F^{n*} \in \mathcal{B}_A \subsetneq \mathcal{S}_A$ and $F \in \mathcal{L}_A$, by \cite[Th. 3.2(i)]{konstantinides:passalidis:2024h} we obtain $F \in \mathcal{S}_A$. Hence, from \cite[Cor. 4.10]{samorodnitsky:sun:2016} we find $F^{n*}(x\,A) \sim n\,F(x\,A)$, therefore, from Proposition \ref{pr.P.2.1}, is implied that $F \in \mathcal{B}_A $.
~\halmos

\section{Single big jump principle in randomly stopped sums} \label{sec.P.4}

Now we examine the asymptotic behavior of the probability
\beam \label{eq.P.4.32}
\PP\left[ {\bf S}_{\tau} \in x\,A\right]\,, 
\eeam
for some fixed set  $A \in \mathscr{R}$, with ${\bf S}_{\tau}$ the multivariate randomly stopped sum in relation \eqref{eq.P.1.1} that follows distribution $F_{{\bf S}_{\tau}}$. Let us make clear, that in this section, $\tau$ is a discrete, one-dimensional random variable, non-degenerate to zero, and the $\{{\bf X}^{(i)}\,,\;i \in \bbn \}$ are i.i.d. copies of the vector ${\bf X}$, whose distribution is $F$. Additionally, $\tau$ is independent of $\{{\bf X}^{(i)}\,,\;i \in \bbn \}$. Our goal is to establish the single big jump principle of randomly stopped sum from \eqref{eq.P.1.1}, namely relation \eqref{eq.P.1.2}, under weaker conditions for the variable $\tau$, in comparison with subexponential case in \cite[Th. 4.2]{konstantinides:passalidis:2024g}. Restricting $F$ in class $\mathcal{S}_A^*$ instead of $\mathcal{S}_A$, we relax the condition 
\beao
\E\left[(1+\vep)^{\tau} \right]< \infty\,,
\eeao 
for some $\vep > 0$, in relation \eqref{eq.P.4.33}. Let us observe, that in case $d=1$, with $A=(1,\,\infty) \in \mathscr{R}$, the asymptotic behavior of \eqref{eq.P.4.32}, in the presence of heavy tails, as well as the distribution class of the randomly stopped sum was studied in several papers, see for example \cite{borovkov:2000}, \cite{aleskeviciene:leipus:siaulys:2008}, \cite{denisov:foss:korshunov:2010}, \cite{leipus:siaulys:2012}, \cite{spindys:siaulys:2020} among others.  

In multidimensional case, the results are much sparse. Next result contains the main result of the section and it is  inspired by \cite[Th. 1(ii)]{denisov:foss:korshunov:2010}. We notice that if $d=1$ and $A = (1,\,\infty)$ it coincides with \cite[Th. 1(ii)]{denisov:foss:korshunov:2010}, in the non-negative case.

\bpr \label{th.P.4.1}
Let $A \in \mathscr{R}$, be some fixed set. If $F \in \mathcal{S}_A^*$ and there exists some $c > \mu_{F_A}$, such that 
\beam \label{eq.P.4.33}
\PP\left[ c\,\tau >x \right] =o\left(\PP\left[{\bf X} \in x\,A \right]\right)\,, 
\eeam
and further $\E[\tau] < \infty$, then relation \eqref{eq.P.1.2} is true, and further it holds $F_{{\bf S}_{\tau}} \in \mathcal{S}_A^*$. 
\epr

\pr~
Since we consider non-negative ${\bf X}$, we obtain $\mu_{F_A}=\E\left[Y_A\right]$. Hence, for any $x>0$ it holds
\beam \label{eq.P.4.34} \notag
\PP\left[ {\bf S}_{\tau} \in x\,A\right] &=& \sum_{n=0}^{\infty} \PP\left[{\bf X}^{(1)}+\cdots +{\bf X}^{(n)}\in x\,A  \right]\,\PP[\tau =n]\\[2mm] \notag
&=& \sum_{n=0}^{\infty} \PP\left[\sup_{{\bf p} \in I_A} {\bf p}^{\top}\,\left({\bf X}^{(1)}+\cdots +{\bf X}^{(n)}\right) > x  \right]\,\PP[\tau =n]\\[2mm] \notag
&\leq& \sum_{n=0}^{\infty} \PP\left[Y_A^{(1)}+\cdots +Y_A^{(n)} > x  \right]\,\PP[\tau =n]= \PP\left[ \sum_{i=0}^{\tau} Y_A^{(i)} >x \right] \\[2mm]
&\sim& \E[\tau]\,\PP\left[  Y_A >x \right] = \E[\tau]\,\PP\left[{\bf X} \in x\,A \right]\,, 
\eeam
where in pre-last step we apply \cite[Th. 1(ii)]{denisov:foss:korshunov:2010}. So, through \eqref{eq.P.4.34} we obtain the upper bound for relation \eqref{eq.P.1.2}.

For the lower bound, for any $n \in \bbn$, because of the fact that ${\bf X}$ is non-negative and the set $A$ is increasing, through the Bonferroni inequality we find that it holds
\beao
\PP\left[ {\bf S}_{n} \in x\,A\right] &\geq&   \PP\left[\bigcup_{i=1}^n \{{\bf X}^{(i)}\in x\,A \} \right]\\[2mm] 
&\geq&  \sum_{i=1}^n \PP\left[{\bf X}^{(i)} \in  x\,A  \right] -   \sum_{1 \leq i<j\leq n}  \PP\left[{\bf X}^{(i)} \in  x\,A\,,\; {\bf X}^{(j)} \in  x\,A \right]\\[2mm] 
&=& n\,\PP\left[Y_A > x  \right] - n\,(n-1)\,\left(\PP\left[Y_A > x  \right] \right)^2  \sim  n\,\PP\left[Y_A > x  \right] =n\,\PP\left[ {\bf X} \in x\,A \right]\,. 
\eeao
Therefore we find
\beam \label{eq.P.4.35} 
\liminf \dfrac{\PP\left[ {\bf S}_{n} \in x\,A\right]}{ n\,\PP\left[ {\bf X} \in x\,A \right]} \geq 1\,, 
\eeam
hence via Fatou's Lemma, since $\E[\tau] < \infty$, it holds
\beam \label{eq.P.4.36} \notag
\PP\left[ {\bf S}_{\tau} \in x\,A\right] &=& \sum_{n=0}^{\infty} \PP\left[ {\bf S}_{n} \in x\,A\right]\,\PP[\tau =n] \\[2mm]
&\gtrsim& \sum_{n=0}^{\infty} n\,\PP\left[ {\bf X} \in x\,A \right]\,\PP[\tau =n]=\E[\tau] \,\PP\left[{\bf X} \in x\,A \right]\,,
\eeam
from where, in combination with \eqref{eq.P.4.34}, is implied  \eqref{eq.P.1.2}. Further from relation \eqref{eq.P.1.2}, the inclusion $F \in \mathcal{S}_A^*$, by Proposition \ref{pr.P.2.1}(ii), we obtain $F_{{\bf S}_{\tau}} \in \mathcal{S}_A^*$.  
~\halmos

\bre \label{rem.P.4.1}
We observe that relation \eqref{eq.P.4.36} holds for any distribution $F$, not only from class $\mathcal{S}^*$, and independently of relation \eqref{eq.P.4.33}. Namely, under the only condition $\E[\tau] < \infty$ we find
\beao 
\liminf \dfrac{\PP\left[ {\bf S}_{\tau} \in x\,A\right]}{ \E[\tau] \,\PP\left[{\bf X} \in x\,A \right]} \geq 1\,.
\eeao
\ere

From the few papers on the topic of multivariate randomly stopped sums, the most focus on the class $MRV$. \cite{das:fasenhartmann:2023} focuses on this topic, with strict enough conditions for the random variable $\tau$, that has been relaxed in \cite[Th. 4.3(iii)]{konstantinides:passalidis:2024g} (but only for sets $A \in \mathscr{R}$).

The next corollary, relaxes further the conditions on $\tau$, through relation \eqref{eq.P.4.33}, however here we consider independent vectors. In the following result we obtain an more explicit relation, than this from \eqref{eq.P.1.2}. 

\bco \label{cor.P.4.1}
Let $A \in \mathscr{R}$, be some fixed set. If $F \in MRV(\alpha,\,V,\,\mu)$, with $\alpha > 1$, and there exists some $c > \alpha$, such that relation \eqref{eq.P.4.33} is valid, and furthermore it holds $\E[\tau] < \infty$, then we find  
\beam \label{eq.P.4.37}
\PP\left[{\bf S}_{\tau} \in x\,A \right] \sim \E[\tau]\,\mu(A)\,\bV(x)\,.
\eeam
\eco

\pr~
Firstly, since $F \in MRV(\alpha,\,V,\,\mu)$, with $\alpha > 1$, we notice that $\mu_{F_A}= \E\left[Y_A\right] < \infty$. Further, from the proof of \cite[Prop. 4.14]{samorodnitsky:sun:2016}, it follows that for all $A \in \mathscr{R}$ it holds $\mu(\partial A)=0$ and $\mu(A) \in (0,\,\infty)$. Thus from Proposition \ref{th.P.4.1} we obtain
\beao
\PP\left[{\bf S}_{\tau} \in x\,A \right] \sim \E[\tau]\,\PP\left[{\bf X} \in x\,A \right] \sim \E[\tau]\,\mu(A)\,\bV(x)\,.
\eeao 
~\halmos

\bre \label{rem.P.4.2}
In the previous corollary, if ${\bf S}_{\tau} ^A:=\sup \left\{u\;:\;{\bf S}_{\tau}  \in u\,A \right\}$ with distribution $F_{S_{\tau}}^A$, then by \eqref{eq.P.4.37} we find that $F_{S_{\tau}}^A \in \mathcal{R}_{-\alpha}$. However, from this is not necessary implied that $F \in MRV(\alpha,\,V,\,\mu)$.
\ere

\section{Vector type precise large deviations principle} \label{sec.P.5}

Here we study asymptotic expressions for the precise large deviations of non-random and random sums in relations \eqref{eq.P.1.3}, that means we focus on asymptotic behavior of the probabilities
\beao
\PP\left[ {\bf S}_{n} \in x\,A\right]\,,\qquad \PP\left[ {\bf S}_{N(t)} \in x\,A\right]\,,
\eeao
as $\nto$, or as $t\to \infty$, respectively. In this section we consider the $\{{\bf X}^{(i)}\,,\;i \in \bbn \}$ as i.i.d. copies of the vector ${\bf X}$, whose distribution is $F$. Additionally in case of random sums we suppose that the counting process $\{N(t)\,,\;t \geq 0\}$ is independent of $\{{\bf X}^{(i)}\,,\;i \in \bbn \}$.

The precise large deviations principle of random sums in one-dimensional set up with $A =(1,\,\infty)$, is a well studied topic, with wide spectrum of appications, see \cite{klueppelberg:mikosch:1997}, \cite{tang:2006b}, \cite{wang:wang:cheng:2006}, \cite{loukissas:2012}, \cite{he:cheng:wang:2013} among others.

In multivariate set up, the most of the existing papers examine either the joint distribution tail or the sum of $d$-lines of business, that means with $A_1$ from relation \eqref{eq.P.2.17}, see for example \cite{wang:wang:2007}, \cite{chen:cheng:2024}, \cite{geng:wang:zhu:2025}, usually is called sum-type precise large deviations principle. For large deviations under $MRV$ and their applications, we refer \cite{hult:lindskog:mikosch:samorodnitsky:2005} and \cite{bhattacharya:palmowksi:zwart:2022}.

In the first main result we find uniform asymptotic estimation for the precise large deviations of non-random sums. It is worth to notice that this result even in one-dimensional subcase with $A=(1,\,\infty)$, generalizes previous publications, see in Remark \ref{rem.P.5.1} below. We note that, the function $h$, is the insensitivity function of distribution $F_A$.

\bth \label{th.P.5.1}
Let some fixed set $A \in \mathscr{R}$. If $F \in \mathcal{S}_{*,A}$, then it holds
\beam \label{eq.P.5.39}
\PP\left[ {\bf S}_n \in x\,A \right] \sim n\,\PP\left[ {\bf X} \in x\,A \right] \,,
\eeam 
as $\nto$, uniformly for $x\geq h^{\leftarrow}\left[n\left(\mu_{F_A}+1\right) \right]$.
\ethe

\pr~
Let us start with estimation of the upper bound in \eqref{eq.P.5.39}. From \cite[Lem. 4.9]{samorodnitsky:sun:2016}, we obtain
\beam \label{eq.P.5.40}
\PP\left[ {\bf S}_n \in x\,A \right] \leq \PP\left[ \sum_{i=1}^n Y_A^{(i)} > x\right] \,,
\eeam
for any $x>0$ and any $n \in \bbn$.

We define a sequence of i.i.d. random variables $\left\{\widetilde{Y}_A^{(i)}\,,\;i \in \bbn \right\}$ with common distribution $\widehat{F}_A$, where 
\beao
\widetilde{Y}_A=Y_A - \mu_{F_A} -1\,,
\eeao 
with mean 
\beao
\E\left[\widetilde{Y}_A\right]=-1\,.
\eeao

We observe that since $F_A \in \mathcal{S}_* \subsetneq \mathcal{L}$, we obtain 
\beao
\overline{\widehat{F}}_A(x) \sim \bF_A(x)\,,
\eeao 
and hence $\widehat{F}_A  \in \mathcal{S}_*\subsetneq\mathcal{S} $, because $F_A  \in \mathcal{S}_* $ and by closure property of class $ \mathcal{S}$ with respect to strong tail equivalence we find $\widehat{F}_A  \in \mathcal{S}$, and from \cite[Lem. 3]{korshunov:2002} it is implied this result. 

Now, from \eqref{eq.P.5.40} in first step and \cite[Th.]{korshunov:2002} in fourth step, we obtain
\beam \label{eq.P.5.41}
&&\PP\left[ {\bf S}_n \in x\,A \right] \leq \PP\left[ \sum_{i=1}^n Y_A^{(i)} > x\right] =\PP\left[ \sum_{i=1}^n \widetilde{Y}_A^{(i)} > x- n\,\mu_{F_A} -n \right]\\[2mm] \notag
&& \leq \PP\left[ \max_{1\leq k \leq n}\sum_{i=1}^k \widetilde{Y}_A^{(i)} > x- n\,\mu_{F_A} -n \right] \sim \dfrac 1{\left|\E\left[\widetilde{Y}_A\right]\right|}\int_{x- n\,\mu_{F_A} -n}^{x- n\,\mu_{F_A} -n+ n\,\left|\E\left[\widetilde{Y}_A\right]\right|}\overline{\widehat{F}}_A(y)\,dy\\[2mm] \notag
&&=\int_{x- n\,\mu_{F_A} -n}^{x- n\,\mu_{F_A} }\bF_A(y+\mu_{F_A}+1)\,dy \leq n\,\bF_A(x- n\,\mu_{F_A} -n+\mu_{F_A}+1)\\[2mm] \notag
&&\leq n\,\bF_A\left(x- n\,\left(\mu_{F_A}+1\right)\right) \leq n\,\bF_A(x- h(x)) \sim n\,\bF_A(x) = n\,\PP[{\bf X} \in x\,A]  \,,
\eeam
as $\xto$, where in the eighth step we use the assumption $x \geq h^{\leftarrow}\left[n\left(\mu_{F_A}+1\right) \right] $. Hence, from \eqref{eq.P.5.41} we conclude
\beam \label{eq.P.5.42}
\limsup_{\nto} \sup_{x \geq  h^{\leftarrow}\left[n\left(\mu_{F_A}+1\right) \right]} \dfrac{\PP\left[ {\bf S}_n \in x\,A \right]}{ n\,\PP[{\bf X} \in x\,A]} \leq 1  \,.
\eeam

Next, we proceed to the lower bound for relation \eqref{eq.P.5.39}. We should mention that it follows easily by \cite[Th. 6.1]{konstantinides:passalidis:2024h}, taking into account that 
\beao
\mu_{F_A}=\E[Y_A] < \infty\,,
\eeao 
and \cite[Rem. 2.2]{he:cheng:wang:2013}. However, for sake of completeness we provide an alternative approach.

Since, ${\bf X}$ is non-negative and the set $A$ is increasing, we obtain through Bonferroni's inequality, that it holds
\beam \label{eq.P.5.43}
&&\PP\left[ {\bf S}_n \in x\,A \right] \geq \PP\left[ \bigcup_{i=1}^n \{ {\bf X}^{(i)} \in x\,A \} \right] \geq  \\[2mm] \notag
&&\sum_{i=1}^n \PP\left[  {\bf X}^{(i)} \in x\,A \right] -  \sum_{1\leq i<j\leq n}  \PP\left[  {\bf X}^{(i)} \in x\,A\,,\; {\bf X}^{(j)} \in x\,A \right] =: n\,\PP\left[  {\bf X} \in x\,A \right]- J_1(x,n)\,.
\eeam
for any $x>0$ and any $n\in \bbn$.

We can see that since $\mu_{F_A}=\E[Y_A] < \infty$, it follows that $x\,\PP\left[ Y_A > x \right] \to 0$, as $\xto$, which further implies that for any $\delta>0$, there exists some sufficiently large $x_0=x_0(\delta)$, such that it holds 
\beao
x\,\PP[Y_A >x] < \delta\,,
\eeao 
for any $x\geq x_0$. Also, because $h(x)=o(x)$, as $\xto$, for any $\epsilon>0$, there exists some sufficiently large $x^*_0=x^*_0(\epsilon)$, such that it holds $h(x)/x < \epsilon$, for any $x\geq x^*_0$. Hence, for any $x\geq x_0 \vee x^*_0$ it holds
\beao
J_1(x,n) &=& \sum_{1\leq i<j\leq n} \PP\left[  {\bf X}^{(i)} \in x\,A\,,\; {\bf X}^{(j)} \in x\,A \right] \\[2mm]
&=& \sum_{j =2}^n  \PP\left[ {\bf X}^{(j)} \in x\,A \right] \sum_{i=1}^{j-1} \PP\left[  {\bf X}^{(i)} \in x\,A\;\big|\; {\bf X}^{(j)} \in x\,A \right] \\[2mm]
&=& \sum_{j =2}^n \dfrac 1{x} \PP\left[ {\bf X}^{(j)} \in x\,A \right] \sum_{i=1}^{j-1} x\,\PP\left[  {\bf X}^{(i)} \in x\,A \right]<\delta \sum_{j =2}^n \dfrac {n-1}{x} \PP\left[ {\bf X}^{(j)} \in x\,A \right] \\[2mm]
&\leq& \delta \sum_{j =2}^n \dfrac {n\,\left(\mu_{F_A}+1\right)}{x}  \PP\left[ {\bf X}^{(j)} \in x\,A \right]  \leq \delta \epsilon \,n\, \PP\left[ {\bf X} \in x\,A \right] \,,
\eeao
where in the pre-last step we use the fact that $\mu_{F_A}>0$, and in last step we take in consideration that $x\geq h^{\leftarrow}\left[n\left(\mu_{F_A}+1\right) \right]$. Hence, from the last relation we can take $\delta \downarrow 0$, $\epsilon \downarrow 0$, to conclude
\beam \label{eq.P.5.44}
\limsup_{\nto} \sup_{x\geq h^{\leftarrow}\left[n\left(\mu_{F_A}+1\right) \right]} \dfrac{J_1(x,n)}{n\,\PP\left[ {\bf X} \in x\,A \right]}=0\,.
\eeam
From relations \eqref{eq.P.5.43}, \eqref{eq.P.5.44} we obtain
\beam \label{eq.P.5.45}
\liminf_{\nto} \inf_{x\geq h^{\leftarrow}\left[n\left(\mu_{F_A}+1\right) \right]} \dfrac{\PP\left[ {\bf S}_n \in x\,A \right]}{n\,\PP\left[ {\bf X} \in x\,A \right]}\geq 1\,.
\eeam
Therefore, from relations \eqref{eq.P.5.42}, \eqref{eq.P.5.45} we get the desired result.
~\halmos

\bre \label{rem.P.5.1}
In the subcase $d=1$, with $A =(1,\,\infty)$ the previous result generalizes \cite[Prop. 2.4]{loukissas:2019}, since it used class $\mathcal{S}_* $ instead of $\mathcal{S}^*$, and any random variables with positive mean instead of $\mu_{F_A} \geq 1$. However, in the previous result we are restricted with respect range of uniformity, over any $x\geq h^{\leftarrow}\left[n\left(\mu_{F_A}+1\right) \right]$ instead over any $x\geq h^{\leftarrow}\left(n\,\mu_{F_A}\right)$. Furthermore, it is easy to see that \cite[Th. 2.2]{loukissas:2019} holds not only on class $\mathcal{S}^*$ but also on $\mathcal{S}_*$
\ere 

Now, we examine precise large deviations for the random sums $S_{N(t)}$, where the counting process $\{N(t)\,,\;t \geq 0\}$ has mean $\lambda(t) =\E[N(t)] < \infty$, for any $t \geq 0$, with $\lambda(t) \to \infty$, as $t\to \infty$. Counting processes $\{N(t)\,,\;t \geq 0\}$ is assumed to be independent  of $\{{\bf X}^{(i)}\,,\;i \in \bbn \}$ and satisfies the following two common conditions, for the one-dimensional precise large deviations, see \cite{klueppelberg:mikosch:1997}, \cite{baltrunas:Leipus:Siaulys:2008} among others.   

\begin{assumption} \label{ass.P.5.1}
It holds
\beao
\dfrac {N(t)}{\lambda(t)} \stackrel{P}{\rightarrow} 1\,,
\eeao
as $t \to \infty$, where $\stackrel{P}{\rightarrow}$ denotes the convergence in probability. 
\end{assumption} 

\begin{assumption} \label{ass.P.5.2}
For any $\delta >0$ let assume that there exists some  $\vep>0$, such that it holds
\beao 
\sum_{n > \left\lfloor (1+\delta)\,\lambda(t) \right\rfloor} (1+ \vep)^n\PP\left[N(t) =n \right] \to 0\,,
\eeao
as $t \to \infty$.
\end{assumption} 

In fact, Assumption \ref{ass.P.5.2} is slightly modified in comparison with the classical assumptions in one-dimensional precise large deviations, where the sum is over $n>(1+\delta)\,\lambda(t)$. This modification is due to interest on precise large deviations with respect to $\left\lfloor \lambda(t) \right\rfloor$ instead of with respect to $\lambda (t)$.

\bth \label{th.P.5.2}
Let some fixed set $A \in \mathscr{R}$. If $F \in \mathcal{S}_{*,A}$, and the counting process $\{N(t)\,,\;t \geq 0\}$ is independent  of $\{{\bf X}^{(i)}\,,\;i \in \bbn \}$ and satisfies Assumptions \ref{ass.P.5.1} and \ref{ass.P.5.2}, then it holds
\beam \label{eq.P.5.46}
\PP\left[ {\bf S}_{N(t)} \in x\,A \right] \sim \left\lfloor \lambda(t) \right\rfloor\,\PP\left[ {\bf X} \in x\,A \right] \,,
\eeam 
as $t \to \infty$, uniformly for $x\geq h^{\leftarrow}\left[\left\lfloor \lambda(t) \right\rfloor\,\left(\mu_{F_A}+1\right) \right]$. 
\ethe

\pr~
We start with the upper bound of \eqref{eq.P.5.46}. At first, we obtain
\beam \label{eq.P.5.47} \notag
\PP\left[ {\bf S}_{N(t)} \in x\,A \right] &=&\left(\sum_{n=1}^{\left\lfloor \lambda(t)\,\left(1+\delta \right) \right\rfloor} + \sum_{n > \left\lfloor \lambda(t)\,\left(1+\delta \right) \right\rfloor} \right)\PP\left[ {\bf S}_{n} \in x\,A \right]\,\PP[N(t)=n]\\[2mm]
&=:& I_1(x,\,t,\,n) + I_2(x,\,t,\,n)\,,
\eeam
for any $x>0$, $t>0$, $\delta >0$.

For the estimation of $I_2(x,\,t,\,n)$, we apply \cite[Prop. 4.12(c)]{samorodnitsky:sun:2016}, to find that for any $\vep = \vep(\delta)>0$ there exists some constant $K>0$, such that it holds
\beam \label{eq.P.5.48} \notag
I_2(x,\,t,\,n) &=& \sum_{n > \left\lfloor \lambda(t)\,\left(1+\delta \right) \right\rfloor} \,\PP\left[ {\bf S}_{n} \in x\,A \right]\,\PP[N(t)=n]\\[2mm] &\leq& K\,\PP\left[ {\bf X} \in x\,A \right] \sum_{n > \left\lfloor \lambda(t)\,\left(1+\delta \right) \right\rfloor} (1+ \vep)^n\,\PP[N(t)=n]\\[2mm]\notag
&=& o\left( \left\lfloor \lambda(t) \right\rfloor\,\PP\left[ {\bf X} \in x\,A \right] \right) \,,
\eeam
as  $t \to \infty$, for any $x>0$, where in the last step we used Assumption \ref{ass.P.5.2}.

For the estimation of $I_1(x,\,t,\,n)$, we apply Theorem \ref{th.P.5.1}, to find that it holds
\beam \label{eq.P.5.49} \notag
I_1(x,\,t,\,n) &\leq& \PP\left[ \sum_{i=1}^{\left\lfloor \lambda(t)\,\left(1+\delta \right) \right\rfloor} {\bf X}^{(i)} \in x\,A \right]\,\sum_{n=1}^{\left\lfloor \lambda(t)\,\left(1+\delta \right) \right\rfloor}  \PP[N(t)=n] \leq \PP\left[ \sum_{i=1}^{\left\lfloor \lambda(t)\,\left(1+\delta \right) \right\rfloor} {\bf X}^{(i)} \in x\,A \right] \\[2mm]
&\sim & \left\lfloor \lambda(t)\,\left(1+\delta \right) \right\rfloor\,\PP\left[ {\bf X} \in x\,A \right] \,,
\eeam
as  $t \to \infty$, uniformly for $x \geq  h^{\leftarrow}\left[\left\lfloor \lambda(t)(1+\delta) \right\rfloor\,\left(\mu_{F_A}+1\right) \right]$. Hence, by \eqref{eq.P.5.49}, letting $\delta \downarrow 0$, we obtain
\beam \label{eq.P.5.50a} 
\limsup_{t \to \infty} \sup_{x \geq  h^{\leftarrow}\left[\left\lfloor \lambda(t) \right\rfloor\,\left(\mu_{F_A}+1\right) \right]} \dfrac {I_1(x,\,t,\,n)}{  \left\lfloor \lambda(t) \right\rfloor\,\PP\left[ {\bf X} \in x\,A \right]} \leq 1 \,.
\eeam
So, by relations \eqref{eq.P.5.48}, \eqref{eq.P.5.50a}, together with \eqref{eq.P.5.47}, we reach to upper asymptotic bound for \eqref{eq.P.5.46}.

For the lower asymptotic bound for \eqref{eq.P.5.46}, it follows directly from \cite[Th. 6.2]{konstantinides:passalidis:2024h}, in combination with condition $\mu_{F_A} < \infty$, and \cite[Rem. 2.2]{he:cheng:wang:2013}. However, for sake of completeness we present an alternative approach. 

Let $\delta > 0$. Then from Assumption \ref{ass.P.5.1}, for any $\vep >0$ we can find some $t^*>0$ sufficiently large, such that it holds
\beao
\PP\left[ \left| \dfrac{N(t)}{\lambda(t)} - 1 \right| < \delta \right] \geq 1- \vep\,,
\eeao 
for any $t \geq t^*$. Hence, we obtain
\beam \label{eq.P.5.51} \notag
\PP\left[ {\bf S}_{N(t)} \in x\,A \right]  &\geq& \sum_{n=\left\lfloor \lambda(t)\,\left(1-\delta \right) \right\rfloor}^{\left\lfloor \lambda(t)\,\left(1+\delta \right) \right\rfloor} \,\PP\left[ {\bf S}_{n} \in x\,A \right]\, \PP[N(t)=n] \\[2mm]\notag
&\geq &  \PP\left[ {\bf S}_{\left\lfloor \lambda(t)\,\left(1-\delta \right) \right\rfloor} \in x\,A \right]\,\PP\left[ \left| \dfrac{N(t)}{\lambda(t)} - 1 \right| < \delta \right]\\[2mm]\notag &\geq& (1-\vep)\,\PP\left[ {\bf S}_{\left\lfloor \lambda(t)\,\left(1-\delta \right) \right\rfloor} \in x\,A \right]\,,
\eeam
for any $x>0$ and any $t\geq t^*$. Consequently, by relation \eqref{eq.P.5.51}, through Theorem \ref{th.P.5.1}, we find
\beao
\PP\left[ {\bf S}_{N(t)} \in x\,A \right] \gtrsim (1-\vep)\,\left\lfloor \lambda(t)\,\left(1-\delta \right) \right\rfloor\,\PP\left[ {\bf X} \in x\,A \right]\,,
\eeao
as $t \to \infty$, uniformly for $x \geq  h^{\leftarrow}\left[\left\lfloor \lambda(t)(1-\delta) \right\rfloor\,\left(\mu_{F_A}+1\right) \right]$. This way, from the arbitrary choice of $\vep>0$ and $\delta > 0$,  letting $\vep \downarrow 0$ and $\delta \downarrow 0$, we obtain
\beao
\liminf_{t \to \infty} \inf_{x \geq  h^{\leftarrow}\left[\left\lfloor \lambda(t) \right\rfloor\,\left(\mu_{F_A}+1\right) \right]} \dfrac{\PP\left[ {\bf S}_{N(t)} \in x\,A \right]}{ \left\lfloor \lambda(t) \right\rfloor\,\PP\left[ {\bf X} \in x\,A \right]} \geq 1\,,
\eeao
that provides the desired lower asymptotic bound of \eqref{eq.P.5.46}.
~\halmos

\bre \label{rem.P.5.2}
If in Theorem \ref{th.P.5.1} we restrict the condition $F \in \mathcal{S}_{*,A}$ into $F \in MRV(\alpha,\,V,\,\mu)$, with $\alpha >1$, then relation \eqref{eq.P.5.39} is reduced to 
\beam \label{eq.P.5.52}
\PP \left[{\bf S}_n \in x\,A \right] \sim n\,\mu(A)\,\bV(x)\,,
\eeam
uniformly for $x\geq h^{\leftarrow}\left(n\,\left(\mu_{F_A}+1\right) \,\right)$, see a similar proof in Corollary \ref{cor.P.4.1}. Correspondingly, in Theorem \ref{th.P.5.2}, if we restrict ourselves into $MRV(\alpha,\,V,\,\mu)$, with $\alpha >1$, then relation \eqref{eq.P.5.46} is reduced to 
\beam \label{eq.P.5.53}
\PP \left[{\bf S}_n \in x\,A \right] \sim \left\lfloor \lambda(t) \right\rfloor\,\mu(A)\,\bV(x)\,,
\eeam
uniformly for $x \geq h^{\leftarrow}\left(\left\lfloor \lambda(t) \right\rfloor\,\left(\mu_{F_A}+1\right) \,\right)$. Relations \eqref{eq.P.5.52}, \eqref{eq.P.5.53}, in comparison with relations \eqref{eq.P.5.39}, \eqref{eq.P.5.46}, have the advantage that the dependence of the components of ${\bf X}$, are described  completely by $\mu$ and the tail decay of ${\bf S}_n$ driven by $V$, which is one of the most important reasons of the domination of $MRV$ in literature. According to our knowledge, the relations \eqref{eq.P.5.52}, \eqref{eq.P.5.53}, in these general forms are not established yet, except only under concrete dependence structures among the components of ${\bf X}$, and for concrete forms of sets $A$. 
\ere

\section{Nonstandard multivariate risk model} \label{sec.P.6}

In this section we consider a nonstandard multivariate risk model, that means the counting process $\{ N(t)\,,\;t \geq 0 \}$ is NOT necessarily renewal. Even though in one-dimensional risk models this set up is still not popular, no matter its theoretical implications or the wider spectrum of practical applications. We mention   the \cite{yang:wang:2012}, for the one-dimensional set up and \cite{chen:wang:wang:2013}, \cite{gao:yang:2014}, for the two-dimensional set up. In these three papers, the distribution class for the claim is the $\mathcal{D} \cap \mathcal{L}$, while they consider a quasi-renewal process, namely weakly dependent
inter-arrival times, that are identically distributed.

Form the other side, the multivariate risk models attracted the attention of several researchers, using mostly i.i.d. claim vectors following distribution $F \in MRV(\alpha,\,V,\,\mu)$ (under the condition $\mu((1,\infty]\times,...,\times(1,\infty]) >0$, which indicates asymptotic dependence and regular variation for the components of claim vector). For renewal risk models we refer to \cite{konstantinides:li:2016}, \cite{li:2016}, \cite{yang:su:2023}, \cite{cheng:konstantinides:wang:2024}, where consider also financial risks.

For distribution classes of claim vectors, larger than $MRV$, we mention \cite{samorodnitsky:sun:2016}, where we find the examination of renewal risk model without interest force, and are provided asymptotic expressions for the infinite-time ruin probability, when the multivariate integrated tail distribution of claim vectors belongs to class $\mathcal{S}_A$. In \cite{konstantinides:passalidis:2024g} was established the asymptotic behavior of the  discounted aggregate claims over a  finite-horizon interval, in the classical risk model with claim distributions from class $\mathcal{S}_A$ and with a c\'{a}dl\'{a}g process for the financial risks.
In the same way in \cite{konstantinides:passalidis:2024j}, with claim distributions from the class $(\mathcal{D}\cap\mathcal{L})_A$ in renewal risk model, with weak dependence among claim vectors, were obtained uniform asymptotic expressions for the ruin probability over a finite-horizon interval.

Next, we establish a local uniform asymptotic expression, for the entrance probability of discounted aggregate claims into some rare-sets, from where follows directly the ruin probability, in a nonstandard risk model  with claim distributions from class $\mathcal{S}_A^*$. We assume that the counting process $\{N(t)\,,\; t\geq 0 \}$, has in some sense lighter distribution tail than that of $F_A$, see Assumption \ref{ass.P.6.2} below.

We start with an analytic representation of the risk model. We suppose that an insurer operates simultaneously $d$-lines of business, with $d \in \bbn$, which share a common counting process of claim-arrivals. The initial capital of the insurer is $x>0$ and it is allocated into the $d$-lines through a deterministic vector ${\bf l}$, with $l_1,\,\ldots,\,l_d >  0$ and
\beao
\sum_{i=1}^d l_i =1\,.
\eeao
We assume that the insurer can invest his surplus into risk free investment, with constant interest force $r\geq 0$. The case $r=0$ corresponds to no investment set up.

Furthermore, we consider that the insurer receives premiums, whose densities are represented by the vector ${\bf p}(t) = (p_1(t),\,\ldots,\,p_d(t))$, for any $t \geq 0$, with $p_i(t)$ denoting the premium density for the $i$-th line of business, for $i=1,\,\ldots,\,d$. For them we adopt the convention $0\leq p_i(t) \leq \Lambda_i$, for any $t\geq 0$ and for some non-negative constants $\Lambda_i$, for $i=1,\,\ldots,\,d$. The $i$-th claim vector ${\bf X}^{(i)} =(X_1^{(i)},\,\ldots,\,X_d^{(i)})$ is non-negative and it can contain zero components (but not all of them). The arrival of ${\bf X}^{(i)}$ happens at time
$\tau_i$, with $i \in \bbn$. 

Let put $\tau_0=0$, then the sequence $\{\tau_i\,,\;i \in \bbn\}$ constitute a counting process $\{N(t)\,,\;t\geq 0\}$, defined as follows
\beam \label{eq.P.6.53} 
N(t) := \sup \{n \in \bbn\;:\;\tau_n \leq t \}\,,
\eeam
for any $t\geq 0$, that by convention has finite mean value
\beam \label{eq.P.6.54} 
\lambda(t) :=\E[N(t) ]= \sum_{i=1}^{\infty} \PP[\tau_i \leq t] < \infty\,,
\eeam
for any $t\geq 0$. As we mentioned, the  counting process $\{N(t)\,,\;t\geq 0\}$ is not necessarily renewal process. More concretely, we consider that the sequence of inter-arrival times $\{\theta_i\,,\;i\in \bbn\}$, where $\theta_i=\tau_i - \tau_{i-1}$, for any $i\in \bbn$, represent sequence of independent but not necessarily identically distributed, positive random variables.

This way, the discounted surplus process for the insurer, ${\bf U}(t)$, at time $t \geq 0$, can be written in the form
\beam \label{eq.P.6.55} 
{\bf U}(t):=\left( 
\begin{array}{c}
U_{1}(t) \\ 
\vdots \\ 
U_{d}(t) 
\end{array} 
\right) =x\,\left( 
\begin{array}{c}
l_{1} \\ 
\vdots \\ 
l_{d}\,
\end{array} 
\right) +\left( 
\begin{array}{c}
\int_{0}^{t} p_{1}(s)e^{-rs}ds \\ 
\vdots \\ 
\int_{0}^{t}p_{d}(s)e^{-rs}ds 
\end{array} 
\right) -\left( 
\begin{array}{c}
\sum_{i=1}^{N(t)}X_{1}^{(i)}e^{-r\tau _{i}} \\ 
\vdots \\ 
\sum_{i=1}^{N(t)}X_{d}^{(i)}e^{-r\tau _{i}} 
\end{array} 
\right)\,. 
\eeam 

Next, we formulate our first Assumption. The independence among the claims, the premiums and the counting process $\{N(t)\,,\;t\geq 0\}$ is quite common either in one-dimensional or in multidimensional risk model, as the condition of i.i.d. claims as well. However, several times the independence between claims and counting process, or the condition of i.i.d. claims seems very restrictive, see for example in \cite{li:2016}, \cite{konstantinides:passalidis:2024j} for more discussions.

\begin{assumption} \label{ass.P.6.1}
Let $A \in \mathscr{R}$ be a fixed set. We assume that the vectors $\left\{{\bf X}^{(k)}\,,\;k \in \bbn \right\}$ are i.i.d. with common distribution $F \in \mathcal{S}_A^*$. We assume that the $\left\{{\bf X}^{(i)}\,,\;i \in \bbn \right\}$, $\{N(t)\,,\;t\geq 0\}$ and $\{{\bf p}(t) \,,\;t\geq 0\}$ are mutually independent. Additionally, we assume that the sequence of inter-arrival times $\{\theta_i\,,\;i \in \bbn\}$ has independent, positive terms.
\end{assumption} 

The assumption that the  $\{\theta_i\,,\;i \in \bbn\}$  are independent but not necessarily identically distributed, escape by a different way from the standard risk models, than the usual one, where are used quasi renewal processes. Under this condition the process $\{N(t)\,,\;t\geq 0\}$ can become  either renewal (but not quasi renewal, due to independence of the $\{\theta_i\,,\;i \in \bbn\}$), or a non-homogeneous renewal process, see for example in \cite{bernackaite:siaulys:2015}.   

Now, we introduce the delayed counting process $\{N^*(t)\,,\;t\geq 0\}$ as follows
\beam \label{eq.P.6.56} 
N^*(t)= \sup \{n \in \bbn\;:\; \tau_n^* \leq t\}\,,
\eeam
for any $t\geq 0$, where
\beam \label{eq.P.6.57} 
\tau_n^*= \sum_{i=2}^{n+1} \theta_i\,,
\eeam
for any $n \in \bbn$. From relation \eqref{eq.P.6.56} in combination with \eqref{eq.P.6.54} we obtain
\beam \label{eq.P.6.58} 
\lambda^*(t)=\E[N^*(t)] =  \sum_{i=1}^{\infty} \PP[\tau_i^* \leq t] < \infty\,,
\eeam 
for any $t\geq 0$. This way we can formulate the following assumption for delayed counting process $\{N^*(t)\,,\;t\geq 0\}$, reminding the mean $\mu_{F_A}= \E{Y_A}$.

\begin{assumption} \label{ass.P.6.2}
Let $\{N^*(t)\,,\;t\geq 0\}$ be a delayed counting process, described through relation \eqref{eq.P.6.56}. We assume that for some constant $c>\mu_{F_A}$ and some  constant $T^*>0$, it holds
\beam \label{eq.P.6.59} 
\sum_{n=0}^{\infty} \dfrac{\PP[N^*(t) \geq n-1]}{\PP[{\bf X} \in c\,n\,A]} < \infty\,,
\eeam 
for any $t \in [0,\,T^*]$.
\end{assumption}

\bre \label{rem.P.6.1}
From the relations \eqref{eq.P.6.56}, \eqref{eq.P.6.57} can be implied that \eqref{eq.P.6.59} is equivalent to 
\beam \label{eq.P.6.60} 
\sum_{n=0}^{\infty} \dfrac{\PP[\theta_2 + \cdots + \theta_n \leq t]}{\PP[{\bf X} \in c\,n\,A]} < \infty\,,
\eeam 
for any  $t \in [0,\,T^*]$. Further, Assumption \ref{ass.P.6.2} implies that the tail of the $\{N(t)\,,\;t\in [0,\,T^*]\}$, is lighter than the tail of $F_A$, in the sense it holds
\beao
\sum_{n=0}^{\infty} \dfrac{\PP[N(t) \geq n]}{\PP[Y_A > c\,n]} < \infty\,.
\eeao
Indeed, since by relation \eqref{eq.P.6.60}  we obtain that it holds
\beao
\infty > \sum_{n=0}^{\infty} \dfrac{\PP[\theta_2 + \cdots + \theta_n \leq t]}{\PP[{\bf X} \in c\,n\,A]} \geq \sum_{n=0}^{\infty} \dfrac{\PP[\theta_1 + \cdots + \theta_n \leq t]}{\PP[{\bf X} \in c\,n\,A]}=\sum_{n=0}^{\infty} \dfrac{\PP[\tau_n \leq t]}{\PP[Y_A > c\,n]}\,,
\eeao
for any $t \in [0,\,T^*]$.
\ere

As in one-dimensional case, the asymptotic behavior of discounted aggregate claims, as for some rare-set $A$, see Remark \ref{rem.P.6.2} below, implies quite easily the ruin probability, as entrance to some ruin-sets. In risk model \eqref{eq.P.6.55}, we denote the last term of the right-hand side, namely the discounted aggregate claims up to moment $t\geq 0$, as 
\beao
{\bf D}_r(t) = \sum_{i=1}^{N(t)} {\bf X}^{(i)} e^{-r\tau_i}\,,
\eeao  
and we are interested for asymptotic estimation of $\PP[{\bf D}_r(t) \in x\,A]$, as $\xto$, for any  $t \leq T$. Namely, the set $x\,A$ is understood as rare-set and the estimations are provided uniformly for $t \in \Lambda_T$,  with $\Lambda_T := \Lambda \cap[0,\,T]$, where $T>0$ some fixed constant and $\Lambda=\{t\;:\;\lambda(t) > 0\}$.

\bre \label{rem.P.6.2}
The set $A \in \mathscr{R}$ can take two main forms, from the aspect to insurance practice. The first form, is described in relation \eqref{eq.P.2.17}, from where it follows that the sum of discounted aggregate claims of the $d$-lines of business, exceed the initial capital. It can be connected with the ruin probability '$\psi_{sum}$'.

The second form can be represented as
\beam \label{eq.P.6.60a}
A_2 = \{{\bf x}\;:\; x_i > b_i\,,\;\exists \; i =1,\,\ldots,\,d \}\,,
\eeam 
with $b_i > 0$, for $i=1,\,\ldots,\,d$, which describes the exceeding of the discounted aggregate claims for only one of the $d$-lines of business. It can be connected with the ruin probability '$\psi_{or}$'. For more information of these ruin probabilities, we refer the reader to \cite{cheng:yu:2019}. Unfortunately, the other two forms of ruin probability, denoted as '$\psi_{and}$' or '$\psi_{sim}$' are not covered by these sets $A$, see in \cite[Rem. 2.2]{konstantinides:passalidis:2024g}. Indeed, for these forms seems better to be approached through  'non-linear' single big jump, see in \cite[Sec. 5]{konstantinides:passalidis:2024a}.   
\ere

\bth \label{th.P.6.1}
Let us consider risk model from \eqref{eq.P.6.55} and we adopt Assumptions \ref{ass.P.6.1} and \ref{ass.P.6.2} for some fixed $T^*$. Then it holds
\beao 
\PP[{\bf D}_r(t) \in x\,A] \sim \int_{0}^{t} \PP\left[{\bf X} \in x\,e^{r\,s}\,A \right]\,\lambda(ds) \,,
\eeao 
uniformly for $t \in \Lambda_T$, with $T \leq T^*$, $T \in \Lambda$.
\ethe

Before the proof of Theorem \ref{th.P.6.1} we need a Lemma, which has merit in its own right.

\ble \label{lem.P.6.1}
Let $A \in \mathscr{R}$ be a fixed set. We assume that the $\{{\bf X}^{(i)}\,,\;i \in \bbn\}$ is a sequence of i.i.d., non-negative random vectors, with common distribution $F \in \mathcal{S}_A$. Then, for any $0<a \leq b < \infty$, it holds 
\beam \label{eq.P.6.65} 
\PP\left[ \sum_{i=1}^n c_i\,{\bf X}^{(i)} \in x\,A\right]\sim \sum_{i=1}^n \PP\left[  c_i\,{\bf X}^{(i)} \in x\,A\right] \,,
\eeam
for any $n\in \bbn$, uniformly for ${\bf c}_n :=(c_1,\,\ldots,\,c_n) \in [a,\,b]^n$.    
\ele

\pr~
For the upper bound, at first we employ \cite[Prop. 2.4]{konstantinides:passalidis:2024g}, and next through \cite[Prop. 5.1]{tang:tsitsiashvili:2003} (or through \cite[Lem. 3.1]{hao:tang:2008}), we obtain
\beam \label{eq.P.6.66} 
\PP\left[ \sum_{i=1}^n c_i\,{\bf X}^{(i)} \in x\,A\right] &\leq& \PP\left[ \sum_{i=1}^n c_i Y_A^{(i)} > x \right] \\[2mm] \notag
&\sim& \sum_{i=1}^n\PP\left[  c_i Y_A^{(i)} > x \right]=\sum_{i=1}^n\PP\left[ c_i {\bf X}^{(i)} \in x\,A\right] ,
\eeam
uniformly for ${\bf c}_n  \in [a,\,b]^n$.

For the lower bound of \eqref{eq.P.6.65}, since the ${\bf X}^{(i)}$ are non-negative,  and the set $A$ is increasing, for any $x>0$, we find 
\beam \label{eq.P.6.67}
\PP\left[ \sum_{i=1}^n c_i\,{\bf X}^{(i)} \in x\,A\right]&\geq& \PP\left[ \bigcup_{i=1}^n \{c_i\,{\bf X}^{(i)} \in x\,A\} \right]  \\[2mm] \notag
&\geq& \sum_{i=1}^n\PP\left[ c_i\,{\bf X}^{(i)} \in x\,A\right]- \sum_{1\leq i<j\leq n} \PP\left[ c_i\,{\bf X}^{(i)} \in x\,A\right]\,\PP\left[ c_j\,{\bf X}^{(j)} \in x\,A\right]  \\[2mm] \notag
&\sim&  \sum_{i=1}^n\PP\left[ c_i\,{\bf X}^{(i)} \in x\,A\right]  \,,
\eeam
for any  ${\bf c}_n  \in [a,\,b]^n$. Hence, by relations \eqref{eq.P.6.66} and \eqref{eq.P.6.67} we get \eqref{eq.P.6.65}, uniformly for ${\bf c}_n  \in [a,\,b]^n$.
~\halmos

Now, we are ready for the proof of Theorem \ref{th.P.6.1}.

{\bf Proof of Theorem  \ref{th.P.6.1}}~
At first, from Assumption \ref{ass.P.6.1}, for some integer $M \in \bbn$ it holds
\beam \label{eq.P.6.68}
&&\PP[{\bf D}_r(t) \in x\,A] =\PP\left[ \sum_{i=1}^{N(t)} {\bf X}^{(i)}\,e^{-r\tau_i} \in x\,A\right]\\[2mm] \notag
&&= \left( \sum_{n=1}^{M} +  \sum_{n=M+1}^{\infty} \right) \PP\left[  \sum_{i=1}^{n} {\bf X}^{(i)}\,e^{-r\tau_i} \in x\,A \,,\; N(t)=n \right] := J_1(x,\,t,\,M) + J_2(x,\,t,\,M) \,,
\eeam
for any $x>0$ and any $t \in \Lambda_T$.

Let us start with estimation of $J_2(x,\,t,\,M)$. For some constant $c>\mu_{F_A}$, by Proposition \ref{pr.P.2.2} we obtain
\beao
J_2(x,\,t,\,M) &\leq&  \sum_{n=M+1}^{\infty} \PP\left[  \sum_{i=1}^{n} {\bf X}^{(i)}\,e^{-r\tau_1} \in x\,A \,,\; \tau_n\leq t \right] \\[2mm]
&=&  \sum_{n=M+1}^{\infty} \int_0^t \PP\left[  \sum_{i=1}^{n} {\bf X}^{(i)}\,e^{-r\,s} \in x\,A\right]\, \PP\left[ \theta_2 +\cdots + \theta_n \leq t -s \right]\,\PP\left[ \tau_1 \in ds \right] \\[2mm]
&\leq&  K\,\sum_{n=M+1}^{\infty} \int_0^t \dfrac{ \PP\left[ \theta_2 +\cdots + \theta_n \leq t -s \right]}{\PP[{\bf X} \in c\,n\,A]}\PP\left[   {\bf X}\,e^{-r\,s} \in x\,A\right]\,\PP\left[ \tau_1 \in ds \right] \\[2mm]
&\leq&  K\,\left(\sum_{n=M+1}^{\infty} \dfrac{ \PP\left[ \theta_2 +\cdots + \theta_n \leq t \right]\,}{\PP[{\bf X} \in c\,n\,A]}\right) \int_0^t \PP\left[   {\bf X} \in x\,e^{r\,s}\,A\right]\,\lambda(ds)\,,
\eeao
for all $t \in \Lambda_T$, where in the second step we take into consideration that the sequence $\{\theta_i,\, i\in\mathbb{N}\}$ has independent terms. Next, from Assumption \ref{ass.P.6.2}, let remind also \eqref{eq.P.6.60}, for any $\vep >0$, there exists some $M>0$ large enough, such that it holds
\beam \label{eq.P.6.69}
J_2(x,\,t,\,M) \leq K\,\vep \int_0^t \PP\left[   {\bf X} \in x\,e^{r\,s}\,A\right]\,\lambda(ds)\,,
\eeam
for all $t \in \Lambda_T$.

Let us continue with estimation of $J_1(x,\,t,\,M)$. By Lemma \ref{lem.P.6.1} we obtain 
\beao
&&J_1(x,\,t,\,M) = \sum_{n=1}^{M} \PP\left[  \sum_{i=1}^{n} {\bf X}^{(i)}\,e^{-r\tau_i} \in x\,A \,,\; N(t) =n \right] \\[2mm]
&& = \sum_{n=1}^{M} \int...\int_{\{0\leq \tau_1 \leq \cdots \leq \tau_n \leq t  < \tau_{n+1}\}} \PP\left[  \sum_{i=1}^{n} {\bf X}^{(i)}\,e^{-r\,s_i} \in x\,A\right]\, \PP\left[ \tau_1 \in ds_1,\,\ldots,\, \tau_{n+1} \in ds_{n+1} \right] \\[2mm]
&& \sim \sum_{n=1}^{M} \int...\int_{\{0\leq \tau_1 \leq \cdots \leq \tau_n \leq t  < \tau_{n+1}\}} \sum_{i=1}^{n} \PP\left[ {\bf X}^{(i)}\,e^{-r\,s_i} \in x\,A\right]\, \PP\left[ \tau_1 \in ds_1,\,\ldots,\, \tau_{n+1} \in ds_{n+1} \right] \\[2mm]
&& =\sum_{n=1}^{M}  \sum_{i=1}^{n} \PP\left[ {\bf X}^{(i)}\,e^{-r\,\tau_i} \in x\,A\,,\;N(t) =n \right]\,,
\eeao
uniformly for $t\in \Lambda_T $, where in the third step, we use the inclusion $\mathcal{S}_A^* \subsetneq \mathcal{S}_A$ and $e^{-r\,T}\leq e^{-r s_i} \leq 1$, for any $i=1,\,\ldots,\,n$. Hence, it  holds
\beam \label{eq.P.6.70} \notag
J_1(x,\,t,\,M) &\sim& \left(\sum_{n=1}^{\infty} - \sum_{n=M+1}^{\infty} \right)\, \sum_{i=1}^{n} \PP\left[ {\bf X}^{(i)}\,e^{-r\,\tau_i} \in x\,A\,,\;N(t) =n \right]\\[2mm]
&=&J_{11}(x,\,t,\,M)  - J_{12}(x,\,t,\,M) \,,
\eeam
uniformly for $t\in \Lambda_T $.

For the first term we interchange the order of summation to find
\beam \label{eq.P.6.70a} 
J_{11}(x,\,t,\,M) \sim \sum_{i=1}^{\infty} \PP\left[ {\bf X}^{(i)}\,e^{-r\,\tau_i} \in x\,A\,,\;\tau_{i} \leq t \right]= \int_0^t \PP\left[   {\bf X} \in x\,e^{r\,s}\,A\right]\,\lambda(ds) \,,
\eeam
uniformly, for $t \in \Lambda_T$. For the second term we obtain for all $t \in \Lambda_T$
\beao
J_{12}(x,\,t,\,M) &\leq& \sum_{n=M+1}^{\infty} \sum_{i=1}^n \PP\left[ {\bf X}^{(i)}\,e^{-r\,\tau_1} \in x\,A\,,\;\tau_n \leq t \right]\\[2mm]
&=&  \sum_{n=M+1}^{\infty} \sum_{i=1}^n \int_0^t \PP\left[   {\bf X} \in x\,e^{r\,s}\,A\right]\, \PP\left[ \theta_2 +\cdots + \theta_n \leq t -s \right]\,\PP\left[ \tau_1 \in ds \right] \\[2mm]
&=&  \sum_{n=M+1}^{\infty} n\, \PP\left[ \theta_2 +\cdots + \theta_n \leq t \right]\, \int_0^t \PP\left[   {\bf X} \in x\,e^{r\,s}\,A\right]\,\PP\left[ \tau_1 \in ds \right] \\[2mm]
&\leq&  \sum_{n=M+1}^{\infty} n\, \PP\left[ N^*(t) \geq n-1 \right]\, \int_0^t \PP\left[   {\bf X} \in x\,e^{r\,s}\,A\right]\,\lambda \left( ds \right)  \,.
\eeao
So, because of \eqref{eq.P.6.58}, for any $\vep>0$ we can find some $M>0$ large enough, such that it holds
\beam \label{eq.P.6.72} 
J_{12}(x,\,t,\,M) \leq \vep \int_0^t \PP\left[   {\bf X} \in x\,e^{r\,s}\,A\right]\,\lambda(ds) \,,
\eeam
for all $t\in \Lambda_T $. Therefore, from relations \eqref{eq.P.6.70}, \eqref{eq.P.6.70a}, and \eqref{eq.P.6.72} we find the two-sided asymptotic inequalities
\beam \label{eq.P.6.73} 
(1-\vep)\, \int_0^t \PP\left[   {\bf X} \in x\,e^{r\,s}\,A\right]\,\lambda(ds)\lesssim J_{1}(x,\,t,\,M) \lesssim \int_0^t \PP\left[   {\bf X} \in x\,e^{r\,s}\,A\right]\,\lambda(ds) \,,
\eeam
uniformly for $t\in \Lambda_T $. Substituting by \eqref{eq.P.6.69} and \eqref{eq.P.6.73} into \eqref{eq.P.6.68}, taking into account the arbitrary choice of $\vep>0$, we conclude the proof.
~\halmos

Next, we find out that the asymptotic behavior of the discounted aggregate claims, can be connected directly with the ruin probability. In model \eqref{eq.P.6.55}, as ruin probability we understand  the probability of surplus entrance into some ruin set $L$, which is related to $A$. Following the paths by \cite{hult:lindskog:2006} and by \cite{samorodnitsky:sun:2016}, we provide the following assumption for the ruin sets.

\begin{assumption} \label{ass.P.6.3}
Let $L$ be some ruin set, that is open and decreasing (that means the $-L$ is increasing set) such that $L^{c}$ is convex and ${\bf 0} \in \partial L$. Let us assume that it holds $\lambda\,L =L$ for any $\lambda >0$.
\end{assumption}

\bre \label{rem.P.6.3}
By  \cite[Sec. 5]{samorodnitsky:sun:2016}, we find two of the forms of ruin sets, that comply with Assumption \ref{ass.P.6.3}
\beao
L_1=\left\{ {\bf x}\;:\;\sum_{i=1}^d x_i<0 \right\}\,, \qquad L_2=\left\{ {\bf x}\;:\;x_i<0\,,\;\exists\;i=1,\,\ldots,\,d\right\}\,,
\eeao
that represent the case, when only one line, from the $d$ ones, has negative surplus, and the case, when the summation of all the surpluses from the $d$-lines of business is negative. Further, we know that $A=({\bf l} - L)$, with ${\bf l}=(l_1,\,\ldots,\,l_d)$, and we can see that this form of ruin probability is satisfied for the $A_2$ of \eqref{eq.P.6.60a} with ruin set $L_2$, while it is satisfied for the set $A_1$ in relation \eqref{eq.P.2.17} for the ruin set $L_1$.

In uni-variate case, the usual ruin set $L=(-\infty,\,0)$ satisfies Assumption  \ref{ass.P.6.3} and then $A = (1,\,\infty)= 1 - L$.
\ere

The finite-horizon ruin probability with respect to set $L$, is defined as follows
\beam \label{eq.P.6.74} 
\psi_{r,L}(x,\,t)=\PP\left[{\bf U}(s) \in L\,,\; \exists\;s \in (0,\,t] \right]\,.
\eeam

Thus, if ${\bf U}(t) $ satisfies the representation \eqref{eq.P.6.55}, and due to Assumption \ref{ass.P.6.3} for the set $L$, this probability \eqref{eq.P.6.74} has the representation
\beao
\psi_{r,L}(x,\,t)&=&\PP\left[{\bf D}_r(s) - \int_0^s {\bf p}(y)\,e^{-r y} \,dy\in x\,({\bf l}-L)\,,\; \exists\;s \in (0,\,t] \right]\\[2mm]
&=& \PP\left[{\bf D}_r(s) - \int_0^s {\bf p}(y)\,e^{-r y} \,dy\in x\,A\,,\; \exists\;s \in (0,\,t] \right]\,,
\eeao
for any $t > 0$, where 
\beao
\int_0^t e^{-r\,y}{\bf p}(y)dy :=
\left( 
\begin{array}{c}
	\int_{0-}^{t} e^{-r\,y}\,p_1(y)\,dy \\ 
	\vdots \\ 
	\int_{0-}^{t} e^{-r\,y}\,p_d(y)\,dy 
\end{array} 
\right) \,,
\eeao

\bco \label{cor.P.6.1}
Let $A =({\bf l}-L) \in \mathscr{R}$ a fixed set, with $L$ to satisfy Assumption \ref{ass.P.6.3}. Under the conditions of Theorem \ref{th.P.6.1}, it holds
\beam \label{eq.P.6.76} 
\psi_{r,L}(x,\,t) \sim\int_0^t \PP\left[   {\bf X} \in x\,e^{r\,s}\,A\right]\,\lambda(ds)\,,
\eeam
uniformly for $t\in \Lambda_T$, with $T \leq T^*$, $T \in \Lambda$.
\eco

\pr~
From the assumptions we find that it holds
\beam \label{eq.P.6.36} 
0 \leq \int_{0}^{t} p_{i}(y)\,e^{-ry}\,dy \leq \Lambda_i\,T < \infty\,,
\eeam
for any  $i=1,\,\ldots,\,n$, and $t\in \Lambda_T$, which permits to apply \cite[Lem. 4.3(d)]{samorodnitsky:sun:2016}. Hence, for the lower asymptotic bound of \eqref{eq.P.6.76}, via \eqref{eq.P.6.36} we obtain that it holds
\beao
\psi_{r,L}(x,\,t) &\geq& \PP\left[{\bf D}_r(t) - \int_0^t {\bf p}(y)\,e^{-r y} \,dy\in x\,A \right] \geq \PP\left[{\bf D}_r(t) \in (x+u)\,A \right] \\[2mm] \notag
&\sim& \int_0^t \PP\left[ {\bf X} \in (x+u)\,e^{r\,s}\,A\right]\,\lambda(ds) \sim \int_0^t \PP\left[ {\bf X} \in x\,e^{r\,s}\,A\right]\,\lambda(ds)\,,
\eeao
uniformly for $t\in \Lambda_T $, where 'u' is a positive constant (which froms by \cite[Lem. 4.3(d)]{samorodnitsky:sun:2016}) and in the third step we apply Theorem \ref{th.P.6.1} and in the last step we use the inclusion $F \in \mathcal{S}_A^* \subsetneq  \mathcal{L}_A $.

From the other-hand side from Theorem  \ref{th.P.6.1} we obtain
\beao
\psi_{r,L}(x,\,t) &\leq& \PP\left[{\bf D}_r(t) \in x\,A \right] \sim \int_0^t \PP\left[ {\bf X} \in x\,e^{r\,s}\,A\right]\,\lambda(ds)\,,
\eeao
uniformly for $t\in \Lambda_T$.
~\halmos

\bre \label{rem.P.6.4}
It seems that Theorem \ref{th.P.6.1} and Corollary \ref{cor.P.6.1} are new even in the one-dimensional set up, for $A=(1,\,\infty)$, ${\bf l} =  1$, $L=(-\infty,\,0)$. However, these results express usual asymptotic expressions, used in one-dimensional risk models. The originality stems from the fact that we require 'lighter' tail in the counting process in comparison with the distribution tail of the claims, and the same time with the help of class $\mathcal{S}^*$, we avoid the restriction in the renewal processes. In some sense, the use of alternative multivariate analogue of Kesten inequality, provided by Proposition \ref{pr.P.2.2}, makes class $\mathcal{S}_A^*$ necessary, and can not be extended to class $\mathcal{S}_A$.
\ere

We close this section with a corollary, that provides a more analytic form in Theorem \ref{th.P.6.1} and in Corollary \ref{cor.P.6.1} in $MRV$ case.

\bco \label{cor.P.6.2}
Let $A =({\bf l}-L) \in \mathscr{R}$ be a fixed set, with $L$ to satisfy Assumption \ref{ass.P.6.3}. Under the conditions of Corollary \ref{cor.P.6.1}, with the restriction that $F \in MRV(\alpha,\,V,\,\mu)$, with $\alpha>1$, it holds
\beam \label{eq.P.6.77} 
\PP\left[{\bf D}_r(t) \in x\,A \right] \sim \psi_{r,L}(x,\,t) \sim \mu(A)\,\bV(x)\,\int_0^t e^{- \alpha\,r\,s}\,\lambda(ds)\,,
\eeam
uniformly for $t\in \Lambda_T$, with $T \leq T^*$, $T \in \Lambda$. 
\eco

\pr~
From Theorem \ref{th.P.6.1} and Corollary \ref{cor.P.6.1}, it is enough to show that
\beao
\int_0^t \PP\left[ {\bf X} \in x\,e^{r\,s}\,A\right]\,\lambda(ds) \sim\mu(A)\,\bV(x)\,\int_0^t e^{- \alpha\,r\,s}\,\lambda(ds)\,, 
\eeao
holds uniformly for $t\in \Lambda_T$. Indeed, since $F \in MRV(\alpha,\,V,\,\mu)$, it follows that $F_A \in \mathcal{R}_{-\alpha}$, for any $A \in \mathscr{R}$, due to $\bF_A(x) \sim \mu(A)\,\bV(x)$, with $V \in \mathcal{R}_{-\alpha}$, and $\mu(A) \in (0,\,\infty)$ for any $A \in \mathscr{R}$ (see the proof of \cite[Prop. 4.14]{samorodnitsky:sun:2016}). Hence, by the uniformity of the regular variation, see \cite[Th. 1.5.2]{bingham:goldie:teugels:1987}, and via the inequalities $1 \leq e^{r\,s} \leq e^{r\,T}$, we find that it holds
\beao
\int_0^t \PP\left[ {\bf X} \in x\,e^{r\,s}\,A\right]\,\lambda(ds) &=&  \int_0^t \bF_A( x\,e^{r\,s})\,\lambda(ds) \sim \bF_A(x) \int_0^t e^{- \alpha\,r\,s}\,\lambda(ds) \\[2mm]
&\sim& \mu(A)\,\bV(x)\,\int_0^t e^{- \alpha\,r\,s}\,\lambda(ds)\,,
\eeao 
this conclude the proof.
~\halmos

\noindent \textbf{Acknowledgments.} 
I would like to thank Dr. Hui Xu, for his comments, which substantially improved the text.


\begin{thebibliography}{99}



\bibitem{aleskeviciene:leipus:siaulys:2008}
{\sc Ale\v{s}kevi\v{c}ien\.{e}, A., Leipus, R., \v{S}iaulys, J}\ (2008) 
Tail behavior of random sums under consistent variation with applications to the compound renewal risk model. 
{\em Extremes}, \textbf{11}, 261-279.




\bibitem{asmussen:1998}
{\sc Asmussen, S.}\ (1998)
Subexponential asymptotics for stochastic processes: extremal behaviour, stationary distributions and first passage probabilities.
{\em Ann. Appl. Probab.}, \textbf{8}, 354--374.




\bibitem{athreya:ney:1972}
{\sc Athreya, K.B., Ney, P.E.}\ (1972)
{\em Branching Processes} 
Springer, New York.


\bibitem{baltrunas:Leipus:Siaulys:2008}
{\sc Baltrunas, A., Leipus, R., \v{S}iaulys, J.}\ (2008)
Precise large deviation results for the total claim amount under subexponential claim sizes. 
{\em Statist. Probab. Lett.}, \textbf{78}, 1206--1214.






\bibitem{bernackaite:siaulys:2015}
{\sc Bernackaite, E., \v{S}iaulys, J.}\ (2015)
The exponential moment tail of inhomogeneous renewal process.
{\em Stat. Probab. Lett.}, \textbf{97}, 9--15. 

\bibitem{bhattacharya:palmowksi:zwart:2022}
{\sc Bhattacharya, A., Palmowski, Z., Zwart, B.}\ (2022) 
Persistence of heavy-tailed sample averages: principle of infinitely many big jumps.
{\em  Electr. J. Probab.} \textbf{27}, no. 50, 1--25.

\bibitem{bingham:goldie:teugels:1987} 
{\sc Bingham. N.H., Goldie, C.M., Teugels, J.L.} \ (1987)
{\em Regular Variation}
Cambridge University Press, Cambridge.

\bibitem{borovkov:2000} 
{\sc Borovkov, A.A.} \ (2000)
Large deviation probabilities for random walks with semiexponential distributions.
{\em Siber. Math. J.}, \textbf{41}, no. 6, 1290--1324.

\bibitem{borovkov:borovkov:2008} 
{\sc Borovkov. A.A., Borovkov, K.A.} \ (2008)
{\em Asymptotic Analysis of Random Walks: Heavy-Tailed Distributions}
Cambridge University Press, Cambridge.


\bibitem{buraczewski:damek:mikosch:2016}
{\sc Buraczewski, D., Damek, E., Mikosch, T.}\ (2016)
{\em Stochastic Models with Power-Law Tails} 
Springer, New York.








\bibitem{chen:cheng:2024}
{\sc Chen, Z., Cheng, D.}\ (2024)
Precise large deviations for non-centralized sums of partial sums and random sums of heavy-tailed END random variables. 
{\em Stat. Probab. Lett.}, \textbf{211}, 110134.






\bibitem{chen:wang:wang:2013}
{\sc Chen, Y., Wang, L., Wang, Y.}\ (2013)
Uniform asymptotics for the finite-time ruin probabilities of two kinds of nonstandard bidimensional risk modes. 
{\em J. Math. Anal. Appl.}, \textbf{401}, no. 1, 114--129.




\bibitem{cheng:konstantinides:wang:2024}
{ \sc Cheng, M., Konstantinides, D.G., Wang, D.}\ (2024)
Multivariate regular varying insurance and financial risks in $d$-dimensional risk model. 
{ \em J. Appl. Probab.}, \textbf{61}, no. 4, 1319 -- 1342.

\bibitem{cheng:yu:2019}
{\sc Cheng, D., Yu, C.}\  (2019)
Uniform asymptotics for the ruin probabilities in a bidimensional renewal risk model with strongly subexponential claims.
{\em Stochastics} \textbf{91}, Vol 1. 643--656.

\bibitem{chistyakov:1964}
{ \sc Chistyakov, V.P.}\ (1964)
A theorem on sums of independent positive random variables and its applications to branching random processes.
{ \em Theory Probab. Appl.}, \textbf{9}, 640--648.

\bibitem{cline:resnick:1992} 
{\sc Cline, D.B.H., Resnick, S.}\ (1992)
Multivariate subexponential distributions.
{\em Stoch. Process. Appl.}, \textbf{42}, no.1, 49--72.



\bibitem{cui:wang:xu:2022} 
{\sc Cui, Z., Wang, Y., Xu, H.}\ (2022)
Some positive conclusions related to the Embrechts - Goldie conjecture.
{\em Sib. Math. J.}, \textbf{63}, 179--192.



\bibitem{das:fasenhartmann:2023}
{\sc Das, B., Fasen-Hartmann, V.}\ (2023)
Aggregating heavy-tailed random vectors: from finite sums to Levy processes.
{\em Preprint, arXiv:2301.10423}.


\bibitem{denisov:foss:korshunov:2004}
{\sc Denisov, D., Foss, S., Korshunov, D.}\ (2004)
Tail asymptotics for the supremum of random walk when the mean is not finite. {\em Queueing Systems}, \textbf{46}, 15--33.

\bibitem{denisov:foss:korshunov:2010}
{\sc Denisov, D., Foss, S., Korshunov, D.}\ (2010)
Asymptotics of randomly stopped sums in the presence of heavy tails.
{\em Bernoulli} \textbf{16}, 971--994.





\bibitem{embrechts:goldie:1980}
{\sc Embrechts, P., Goldie, C. M.}\ (1980)
On closure and factorization properties of subexponential and related distributions. 
{\em J. Austral. Math. Soc. (Ser. A)}, \textbf{29}, 243--256.

\bibitem{embrechts:goldie:veraverbeke:1979}
{\sc Embrechts, P., Goldie, C.M., Veraverbeke, N.}\ (1979)
Subexponentiality and infinite divisibility.
{\em  Z. Wahrscheinlichkeitstheorie Verw. Gebiete}, \textbf{49}, 335--347.



\bibitem{feller:1969}
{\sc Feller, W.}\ (1969)
One-sided\quad analogues\quad of\quad  Karamata's\quad regular\quad  variation.
{\em L' enseignement Math\'{e}matique}, \textbf{15}, 107--121.


\bibitem{foss:korshunov:palmowski:2024} 
{\sc Foss, S., Korshunov, D., Palmowski, Z.} \ (2024)
Maxima over random  time intervals for heavy-tailed compound renewal and L\'{e}vy processes.
{\em Stoch. Process. Appl.}, \textbf{176}, 104422.

\bibitem{foss:korshunov:zachary:2013} 
{\sc Foss, S., Korshunov, D., Zachary, S.} \ (2013)
{\em An Introduction to Heavy-Tailed and Subexponential Distributions.}
Springer, New York, 2nd ed.

\bibitem{foss:zachary:2003}
{\sc Foss, S., Zachary, S.}\ (2003)
The maximum on a random time interval of a random walk with long-tailed increments and negative drift.
{\em Ann. Appl. Probab.},  \textbf{13}, no. 1, 37--53.




\bibitem{gao:yang:2014}
{\sc Gao, Q., Yang, X.}\ (2014)
Asymptotic ruin probabilities in a generalized bidimensional risk model perturbed by diffusion with constant force of interest. 
{\em J. Math. Anal. App.}, \textbf{419}, no. 2, 1193--1213.




\bibitem{geng:wang:zhu:2025}
{\sc Geng, B., Wang, S., Zhu, W.}\  (2025)
Vector-type precise large deviations for a nonstandard multidimensional risk model with some arbitrary dependence structures. 
{\em Acta Math. Hungar.}, https://doi.org/10.1007/s10474-024-01501-3.

\bibitem{goldie:1978}
{\sc Goldie, C.M.}\ (1978)
Subexponential distributions and dominated variation tails 
{\em J. Appl. Probab.}, \textbf{15}, 440--442.




\bibitem{haan:resnick:1981} 
{\sc Haan, L. de, Resnick, S.}\ (1981)
On the observation closet to the origin.
{\em  Stoch. Process. Appl.}, \textbf{11}, no. 3, 301--308.




\bibitem{hao:tang:2008} 
{\sc Hao, X., Tang, Q.}\ (2008)
A uniform asymptotic estimate for discounted aggregate claims with subexponential tails.
{\em  Insur. Math. Econom.}, \textbf{43}, 116--120.

\bibitem{he:cheng:wang:2013}
{\sc He, W., Cheng, D., Wang, Y.}\  (2013)
Asymptotic lower bounds of precise large deviations with nonnegative and dependent random variables. 
{\em Stat. Probab. Lett.}, \textbf{83}, 331--338.



\bibitem{hult:lindskog:2006}
{\sc Hult, H., Lindskog, F.}\ (2006)
Heavy-tailed insurance protfolios: buffer capital and ruin probabilities.
{\em Technical report}.

\bibitem{hult:lindskog:mikosch:samorodnitsky:2005}
{\sc Hult, H., Lindskog, F., Mikosch, T., Samorodnitsky, G.}\ (2005)
Functional large deviations for multivariate regularly varying random walks.
{\em Ann. Appl. Probab.}, \textbf{15}, 2651--2680.







\bibitem{kaas:tang:2003}
{\sc Kaas, R., Tang, Q.}\  (2003)
Note on the tail behavior of random walk maxima with heavy tails and negative drift. 
{\em N. Amer. Act. J.}, \textbf{7}, 57--61.



\bibitem{klueppelberg:1988}
{\sc Kl{\"{u}}ppelberg, C.}\ (1988)  
Subexponential distributions and integrated tails.
{\em J. Appl. Probab.}, \textbf{25}, 132--141.

\bibitem{klueppelberg:mikosch:1997}
{\sc Kl{\"{u}}ppelberg, C., Mikosch, T.}\ (1997)  
Large deviations of heavy-tailed random sums with applications in insurance and finance.
{\em J. Appl. Probab.}, \textbf{34}, 293--308.
	

\bibitem{konstantinides:2018} 
{\sc Konstantinides, D.G.} \ (2018)
{\em Risk Theory. A Heavy Tail Approach.}
World Scientific, New Jersey.

\bibitem{konstantinides:leipus:siaulys:2023} 
{\sc Konstantinides, D.G., Leipus, R., \v{S}iaulys, J.}\ (2023)
On the non-closure under convolution for strong subexponential distributions.  
{\em Nonlin. Anal.: Model. and Contr.}, {\bf 28}, No. 1, 97--115.



\bibitem{konstantinides:li:2016}
{\sc Konstantinides, D.G., Li, J.}\ (2016)
Asymptotic ruin probabilities for a multidimensional renewal risk model with multivariate regularly varying claims.
{\em Insur. Math. and Econom.}, \textbf{69}, 38--44.

\bibitem{konstantinides:liu:passalidis:2026} 
{\sc Konstantinides, D.G., Liu, J., Passalidis, C.D.} \ (2025)
Uniform asymptotics for a multidimensional renewal risk model with multivariate subexponential claims.
{\em Scand. Actuar. J.}, p. 1 -- 21.\\ 
DOI: 10.1080/03461238.2025.2584008. 


\bibitem{konstantinides:mikosch:2005}
{\sc Konstantinides, D.G., Mikosch, T.}\ (2005)
Large Deviations and Ruin Probabilities for Solutions to Stochastic Recurrence Equations with Heavy-tailed Innovations.
{\em Ann. Probab.}, \textbf{33}, 1992--2035.

\bibitem{konstantinides:passalidis:2024a} 
{\sc Konstantinides, D.G., Passalidis, C.D.} \ (2025)
Background risk model in presence of heavy tails under dependence.
{\it Nonlin. Anal. Model. Contr.}, \textbf{30}, n. 5, 982 -- 1010.



\bibitem{konstantinides:passalidis:2024g} 
{\sc Konstantinides, D.G., Passalidis, C.D.} \ (2024)
Random vectors in the presence of a single big jump.
{\em Preprint, arXiv:2410.10292}.

\bibitem{konstantinides:passalidis:2024j} 
{\sc Konstantinides, D.G., Passalidis, C.D.} \ (2025)
Uniform asymptotic estimates for ruin probabilities of a multidimensional risk model with c\'{a}dl\'{a}g returns and multivariate heavy-tailed claims.
{\it Insur. Math. Econom.}, \textbf{125}, 103148.

\bibitem{konstantinides:passalidis:2024c} 
{\sc Konstantinides, D.G., Passalidis, C.D.} \ (2025)
A new approach in two-dimensional heavy-tailed distributions.
{\it Ann. Actuar. Scien.}, \textbf{19}, 317 -- 349.

\bibitem{konstantinides:passalidis:2024h} 
{\sc Konstantinides, D.G., Passalidis, C.D.} \ (2025)
Heavy-tailed random vectors: theory and applications.
{\em Preprint, arXiv:2503.12842}.


\bibitem{korshunov:2002}
{\sc Korshunov, D.}\ (2002) 
Large-deviation probabilities for maxima of sums of independent random variables with negative mean and subexponential distribution.
{\em Theor. Probab. Appl.}, \textbf{46}, 355--366.

\bibitem{korshunov:2018}
{\sc Korshunov, D.}\ (2018)
On subexponential tails for the maxima of negatively driven compound renewal and L\'{e}vy processes.
{\em Stoch. Process. Appl.}, \textbf{128}, 1316--1332.


\bibitem{leipus:siaulys:2012}
{\sc Leipus, R., {\v S}iaulys, J.}\ (2012)
Closure of some heavy-tailed distribution classes under random convolution.
{\em Lithuan. Math. J.}, \textbf{52}, 249--258.

\bibitem{leipus:siaulys:2020}
{\sc Leipus, R.,  {\v S}iaulys, J.}\ (2020)
On a closure property of convolution equivalent class of distributions.
{\em J. Math. Anal. Appl.}, \textbf{490}, no. 124226.



\bibitem{leipus:siaulys:konstantinides:2023}
{\sc Leipus, R., \v{S}iaulys, J., Konstantinides, D.G.}\ (2023)
{\em Closure Properties for Heavy-Tailed and Related Distributions: An Overview.}
Springer Nature, Cham Switzerland.



\bibitem{leslie:1989}
{\sc Leslie, J.R.}\ (1989)
On the non-closure under convolution of the subexponential family. 
{\em J. Appl. Probab.}, \textbf{26}, 58--66.  


\bibitem{li:2016}
{\sc Li, J.}\ (2016)
Uniform asymptotics for a multi-dimensional time-dependent risk model with multivariate regularly varying claims and stochastic return. 
{\em Insur. Math. Econom.}, \textbf{71}, 195--204.




\bibitem{li:2022}
{\sc Li, J.}\ (2022)
Asymptotic results on marginal expected shortfalls for dependent risks.
{\em Insur. Math. Econom.}, \textbf{102}, 310--324.





\bibitem{li:tang:wu:2010} 
{\sc Li, J., Tang, Q., Wu, R.}\ (2010)
Subexponential tails of discounted aggregate claims in a time-dependent renewal risk model. 
{\em Adv. Appl. Probab.}, \textbf{42}, no. 4, 1126-1146.



\bibitem{liu:yi:2025}
{\sc Liu, J., Yi, Y.}\ (2025)
Asymptotics for the conditional higher moment coherent risk measure with weak contagion. 
{\em ASTIN bull.}, \textbf{55}, no. 1, 121--143.

\bibitem{loukissas:2012} 
{\sc Loukissas, F.}\ (2012)
Precise large deviations for long-tailed distributions. 
{\em J. Theor. Probab.}, \textbf{25}, 913--924.

\bibitem{loukissas:2019} 
{\sc Loukissas, F.}\ (2019)
Precise large deviations for strong subexponential distributions and applications on a multi risk model. 
{\em Commun. Statist. Theor. Meth.}, \textbf{48}, no. 20, 5175--5190.




\bibitem{man:tang:2024}
{\sc Man, X., Tang, Q.}\ (2024)
Tail risk driven by investment losses and exogenous shocks. 
{\em Astin Bull.}, \textbf{54}, 712--737.




\bibitem{mikosch:samorodnitsky:2000a} 
{\sc Mikosch, T., Samorodnitsky, G.}\ (2000) 
The supremum of a negative drift random walk with dependent heavy-tailed steps. 
{\em Ann. Appl. Probab.}, \textbf{10}, no. 3, 1025--1064.








\bibitem{omey:2006}
{\sc Omey, E.} \ (2006)
Subexponential distribution functions in $R^{d}$.
{\em J. Math. Sci.}, \textbf{138}, no. 1, 5434--5449.




\bibitem{palmowski:pojer:thonhauser:2025}
{\sc Palmowski, Z., Pojer, S., Thonhauser, S.} \ (2025)
Exact asymptotics or ruin probabilities with linear Hawkes arrivals.
{\em Stoch. Process. Appl.}, \textbf{182}, 104571.

\bibitem{qian:geng:wang:2022}
{\sc Qian, H., Geng, B., Wang, S.} \ (2022)
Tail asymptotics of randomly weighted sums of dependent strong subexponential random variables.
{\em Lith. Math. J.}, \textbf{62}, no. 1, 113--122.


\bibitem{resnick:2007}
{\sc Resnick, S.}\ (2007) 
{\em Heavy-Tail Phenomena. Probabilistic and Statistical Modeling.} 
Springer, New York.






\bibitem{samorodnitsky:sun:2016} 
{\sc Samorodnitsky, G., Sun, J.}\ (2016) 
Multivariate subexponential distributions and their applications. 
{\em Extremes}, \textbf{19}, no. 2, 171--196.







\bibitem{spindys:siaulys:2020}
{\sc Spindys, J., {\v{S}}iaulys, J.}\  (2020) 
Regularly distributed randomly stopped sum, minimum, and maximum.
{\em Nonlin. Anal. Model. Contr.}, \textbf{25}, 509--522.



\bibitem{tang:2006b}
{\sc Tang, Q.}\ (2006)
Insensitivity to negative dependence of asymptotic behavior of precise large deviations.
{\em Electron. J. Probab.}, \textbf{11}, 107--120.




\bibitem{tang:tsitsiashvili:2003}
{\sc Tang, Q., Tsitsiashvili, G.}\ (2003)
Randomly weighted sums of subexponential random variables with application to ruin theory.
{\em Extremes}, \textbf{6}, 171--188.



\bibitem{teugels:1975}
Teugels, J.L., (1975)
The class of subexponential distributions.
\emph{ Ann. Probab.}, {\bf 3}, 1000--1011.




\bibitem{wang:wang:2007}
{\sc Wang, S., Wang, W.}\ (2007)
Precize large deviations for sums of random  variables with consistently varying tails in multi-risk models. 
{\em J. Appl. Probab.}, \textbf{44}, no. 4, 889--900.

\bibitem{wang:wang:cheng:2006}
{\sc Wang, Y., Wang, K., Cheng, D.}\ (2006)
Precize large deviations for sums of negatively associated random  variables with common dominatedly varying tails. 
{\em Acta Math. Sin. (Engl. Ser.)}, \textbf{22}, no.6, 1725--1734.


\bibitem{watanabe:2019}
{\sc Watanabe, T.}\ (2019)
The Wiener condition and the conjectures of Embrechts and Goldie.
{\em Ann. Probab.}, \textbf{47}, no. 3, 1221--1239.

\bibitem{watanabe:yamamuro:2010}
{\sc Watanabe, T., Yamamuro, K.}\  (2010) 
Ratio of the tail of an infinitely divisible distribution on the line to thats of its {L}{\'{e}}vy measure. 
{\em Electr. J. Probab.}, \textbf{15},  44--74.

\bibitem{xu:foss:wang:2015}
{\sc Xu, H., Foss, S., Wang, Y.}\ (2015)
Convolution and convolution-root properties of long-tailed distributions.
{\em Extremes}, \textbf{18}, 605--628.


\bibitem{xu:yu:wang:cheng:2024}
{\sc Xu, H., Yu., C., Wang, Y., Cheng, D.}\ (2024)
Closure under infinitely divisible distribution roots and the Embrechts - Goldie conjecture.
{\em Lithuan. Math. J.}, \textbf{64}, 101--114.







\bibitem{yang:su:2023}
{\sc Yang, Y., Su, Q.}\ (2023)
Asymptotic behavior of ruin probabilities in a multidimensional risk model with investment and multivariate regularly varying claims.
{\em J. Math. Anal. Appl.}, \textbf{525}, 127319.

\bibitem{yang:wang:2012}
{\sc Yang, Y., Wang, K.}\ (2012)
Uniform asymptotics for the finite-time and infinite-time ruin probabilities in a dependent risk model with constant interest rate and heavy-tailed claims.
{\em Lithuan. Math. J.}, \textbf{52}, no.1, 111--121.












\end{thebibliography}
\end{document}